\documentclass[leqno]{amsart}

\usepackage{amscd,amssymb,amsxtra,amsmath}
\usepackage[pdftex]{hyperref}
\usepackage[all]{xy}

\input xy
\xyoption{all}

\DeclareMathAlphabet{\scr}{OMS}{rsfs}{m}{n}
\DeclareMathAlphabet{\eus}{U}{eus}{m}{n}

\newcommand{\mc}[1]{\mathcal{#1}}

\newcommand{\tens}[1]{\boxtimes_{\mc{#1}}}
\newcommand{\score}[0]{\underline{\;\;\;}}
\newcommand{\bitens}[3]{\mc{#1}\tens{#2}\mc{#3}}
\newcommand{\Buni}[2]{B_{\mc{#1}, \mc{#2}}}

\newcommand{\cent}[2]{\mc{Z}_{\mc{#2}}(\mc{#1})}
\newcommand{\arb}[1]{\ar@<-.15ex>[#1]\ar@<.15ex>[#1]\ar@<-.075ex>@{-}[#1]\ar@<.075ex>@{-}[#1]\ar[#1]}
\newcommand{\arbdot}[1]{\ar@<-.15ex>@{:}[#1]\ar@<.15ex>@{:}[#1]\ar@<-.075ex>@{:>}[#1]\ar@<.075ex>@{:>}[#1]\ar@{:>}[#1]}

\newcommand{\ind}[0]{\textrm{ind}}
\newcommand{\End}[0]{\textrm{End}}
\newcommand{\Rep}[0]{\textrm{Rep}}

\newcommand{\Forg}[0]{\textrm{Forg}}

\newcommand{\FPdim}{\textrm{FPdim}}

\newcommand{\Aut}{\textrm{Aut}}

\newcommand{\ev}{\textrm{ev}}
\newcommand{\coev}{\textrm{coev}}

\newcommand{\RQq}{\!\!\!\!\!\!\!\!}

\newtheorem{thm}{Theorem}[section]
\newtheorem{lem}[thm]{Lemma}
\newtheorem{prop}[thm]{Proposition}
\newtheorem{cor}[thm]{Corollary}
\theoremstyle{definition}\newtheorem{defn}[thm]{Definition}
\newtheorem{ex}[thm]{Example}

\newtheorem{remark}[thm]{Remark}

\newtheorem{note}[thm]{Note}
\newtheorem{notation}[thm]{Notation}

\begin{document}

\title{relative centers and tensor products of tensor and braided fusion categories}
\author{Justin Greenough}

\email{jrg8@cisunix.unh.edu}
\address{Department of Mathematics and Statistics\\University of New Hampshire}
\address{Department of Mathematics\\University of Massachusetts Boston}
\date{}
\maketitle
\begin{abstract}
In this paper we study the relative tensor product of module categories over braided fusion categories using, in part, the notion of the relative center of a module category. In particular we investigate the canonical tensor category structure and braiding inherited by the relative tensor product when module categories are themselves tensor/braided. As a basic example we show that the category of representations of the fibre product of finite groups may be expressed in terms of the relative product and show how to multiply braided fusion categories arising from pre-metric groups. Also, we consider relative centers of braided fusion categories and look at the relative tensor product over M\"{u}ger centers. We finish with an in-depth example hinting at a categorification of conjugacy classes for finite groups. 
\end{abstract}
\thispagestyle{empty}
\setcounter{tocdepth}{1}
\tableofcontents
\section{Introduction and main results}\label{IntroductionSection}
The relative tensor product of module categories (defined in Section \ref{PrelimNotats} below) is defined in a way which naturally extends the tensor product of abelian categories as defined by Deligne in \cite{CatsTann}. This construction appears in \cite{ENO:homotop}, and in \cite{JG} where it is studied from the standpoint of 2-categories. This construction is important in describing and classifying various category theoretical notions, e.g. it turns out that de-equivariantization can be described concisely in terms of the relative tensor product (see \textit{loc. cit.} and \cite{braidedI}). 

In this paper we are interested in examining the relative tensor product of module categories which also possess extra structures or properties. In particular we are interested in asking questions along the following lines: If $\mc{D}$-bimodule categories $\mc{C}_1, \mc{C}_2$ also happen to be tensor categories does $\mc{C}_1\tens{D}\mc{C}_2$ possess the structure of a tensor category in some canonical way? If $\mc{C}_1, \mc{C}_2$ are braided, when does the tensor product have a canonical braiding? What is its center? Clearly it is possible to formulate many interesting questions. We hope to answer some of them here. 

The following concept is essential to meaningfully formulating and addressing these questions and is therefore also essential to this paper.
\begin{defn}[\cite{braidedI}]\label{overnessdef} Let $\mc{C}$ be a tensor category and let $\mc{D}$ be a symmetric fusion category. We say $\mc{C}$ is tensor \textit{over} $\mc{D}$ if there is a braided tensor inclusion $\sigma:\mc{D}\hookrightarrow \mc{Z}(\mc{C})$ such that the composition $\mc{D}\hookrightarrow \mc{Z}(\mc{C})\rightarrow \mc{C}$ is faithful. 
\end{defn}
The situation described in Definition \ref{overnessdef} allows us to identify $\mc{D}$ with its image in the center $\mc{Z}(\mc{C})$, and thus as a tensor subcategory of $\mc{C}$ which also has a braiding (its forgetful image). This inclusion gives $\mc{C}$ the structure of a $\mc{D}$-bimodule category. 
In the sequel we will identify the image of any object in $\mc{D}$ under $\sigma$ with its image in $\mc{C}$ without hesitation.

The first result illustrating how structures exhibited by the tensorand categories descend to the relative tensor product is the following theorem, which is proved in Section \ref{sectionRelCentMonCat}. Here $\mc{Z}_\mc{D}(\mc{C})$ is a certain generalization of the categorical center (see Definition \ref{defCent}).
\begin{thm}\label{tensRelCentThm} Let $\mc{C}$ be tensor over $\mc{D}$. Then the  category $\cent{C}{D}$ is tensor over $\mc{D}$. If, in addition, $\mc{C}$ is fusion then $\cent{C}{D}$ is also fusion. 
\end{thm}
Theorem \ref{tensRelCentThm} immediately yields the corresponding result describing the monoidal structure inherited by the relative tensor product, also proved in Section \ref{sectionRelCentMonCat}. 
\begin{cor}\label{monoidalityOfTensor} Let $\mc{C}_1, \mc{C}_2$ be tensor categories each of which is tensor over $\mc{D}$. Then the relative tensor product $\mc{C}_1\tens{D}\mc{C}_2$ has a canonical tensor category structure. If both $\mc{C}_i$ are fusion then the relative tensor product is also fusion.
\end{cor}
Theorem \ref{tensRelCentThm} and Corollary \ref{monoidalityOfTensor} provide new ways to construct fusion categories. In Section \ref{relCentBraidSect} we discuss when the relative center of a braided tensor category, and thus by corollary the relative tensor product of a pair of braided fusion categories, is itself braided. One might hope that, for $\mc{C}$ over $\mc{D}$, the relative center $\cent{C}{D}$ is braided in some natural way. Unfortunately it's not evident under what conditions a braiding exists. However, we are able to show that the relative tensor product is braided if we assume a version of the phenomenon described in Definition \ref{overnessdef} appropriate for braided categories.  
\begin{defn}[\cite{braidedI}]\label{braidingOVERdef} Let $\mc{C}, \mc{D}$ be braided tensor categories. We say $\mc{C}$ is braided \textit{over} $\mc{D}$ if there is a braided inclusion $\mc{D}\hookrightarrow\mc{C}'$.
\end{defn}
Here $\mc{C}'$ is the M\"{u}ger center of $\mc{C}$ (see Definition \ref{centralizerDefn}). In Section \ref{relCentBraidSect} we explicitly describe the braiding of the relative tensor product $\mc{C}_1\tens{D}\mc{C}_2$ where $\mc{C}_i$ are braided over $\mc{D}$, giving a braided version of Corollary \ref{monoidalityOfTensor}.
\begin{thm}\label{braiding of relative center thm} Suppose that $\mc{C}_1, \mc{C}_2$ are braided over $\mc{D}$ as in Definition \ref{braidingOVERdef}. Then $\mc{C}_1\tens{D}\mc{C}_2$ has a canonical braiding coming from the braidings in $\mc{C}_1, \mc{C}_2$ such that the universal balanced functor $B_{1,2}:\mc{C}_1\tens{}\mc{C}_2\rightarrow\mc{C}_1\tens{D}\mc{C}_2$ is braided. 
\end{thm}
In Section \ref{ExampleSection} we give a few examples of the relative product of known braided categories in familiar situations exemplifying Theorem \ref{braiding of relative center thm} and its corollaries. Explicitly we describe the relative tensor product of pointed braided fusion categories (\`{a} la pre-metric groups, discussed in Section \ref{PreMetGrpSec}) and show how the representation category of the fibre product of finite groups occurs naturally as the tensor product of categories of group representations. 

In Section \ref{AlmostNonDegSect} we turn our attention to the centers of braided fusion categories, applying general results developed in \ref{LemmasSection}. For braided fusion category $\mc{C}$ the center $\mc{Z}(\mc{C})$ contains, as a braided fusion subcategory, the relative tensor product $\widetilde{\mc{C}}:=\mc{C}\boxtimes_{\mc{C}'}\mc{C}^{rev}$ over the M\"{u}ger center $\mc{C}'$. As a basic result we show that the center decomposes as a $\widetilde{\mc{C}}$-bimodule category:
\begin{equation}\label{blahblah}
\mc{Z}(\mc{C})\simeq\widetilde{\mc{C}}\oplus\mc{R}
\end{equation}
where $\mc{R}$ is some non-trivial $\widetilde{\mc{C}}$-bimodule subcategory of $\mc{Z}(\mc{C})$. As an example we determine $\mc{R}$ explicitly for $\mc{C}=\Rep(G)$ where $G$ is any finite group. In this case $\mc{R}$ decomposes into a sum of subcategories $\mc{C}_a$ of $\mc{Z}(\Rep(G))$ indexed by a complete set of representatives $R$ of conjugacy classes in $G$. We show that the set $\{\mc{C}_a|a\in R\}$ forms the basis of a $\mathbb{Z}$-based ring (multiplication coming from $\boxtimes_{\Rep(G)}$) isomorphic to the center of the group ring $\mathbb{Z}[G]$. The hope is that, for general braided fusion $\mc{C}$, the decomposition of $\mc{R}$ in (\ref{blahblah}) will provide a categorical interpretation of conjugacy classes. 
\subsection*{Acknowledgements} The author wishes to thank Dmitri Nikshych and Alexei Davydov for valuable discussion and comments, and in particular would like to thank Dmitri Nikshych for spotting errors in earlier versions.
\section{Preliminaries, notation}\label{PrelimNotats} In this section we fix notation and explain (sometimes sketch) ideas essential to the paper. Almost nothing in this section is new and we urge the reader to find more detailed accounts in the references provided.  All categories are assumed abelian, finite, and $k$-linear for $k$ some fixed field (typically of characteristic $0$, though this assumption is not always necessary). All functors are assumed additive and $k$-linear. By \textit{tensor category} we mean a finite $k$-linear abelian rigid monoidal category with the property that $\End(1)\simeq k$, even though for many results neither finiteness nor rigidity is necessary. A tensor category is called \textit{fusion} if it is semisimple with a finite number of simple objects.
\subsection{Module categories} We begin by providing a definition of central importance to this paper. For further details and discussion we recommend the beautiful paper \cite{O:MC}.
\begin{defn} A \textit{right module category} over tensor category $\mc{C}$ is a category $\mc{M}$ together with a bifunctor $\otimes:\mc{M}\times\mc{C}\rightarrow\mc{M}$ and a family of natural isomorphisms $\mu_{M, X, Y}:(M\otimes X)\otimes Y\rightarrow M\otimes (X\otimes Y)$, $r_M:M\otimes 1\rightarrow M$ for $X, Y\in\mc{C}$ and $M\in\mc{M}$ subject to certain natural coherence axioms (one pentagon and one triangle). Similarly one defines the structure of \textit{left} module category. If the structure maps are identity we say $\mc{M}$ is \textit{strict} as a module category over $\mc{C}$.
\end{defn}
\begin{defn}\label{BimoduleCatDef} Let $\mc{C}$, $\mc{D}$ be tensor categories. A category $\mc{M}$ is called a $\mc{C}, \mc{D}$-\textit{bimodule category} if it is a $\mc{C}\tens{}\mc{D}^{rev}$-module category. Here $\mc{D}^{rev}$ signifies the category which is equivalent to $\mc{D}$ as an abelian category but with the tensor structure $X\otimes^{rev}Y:=Y\otimes X$. Any such bimodule category $\mc{M}$ has the obvious left $\mc{C}$-module and right $\mc{D}$-module category structures defined by 
\begin{equation*}
X\otimes M:=(X\tens{}1)\otimes M,\qquad M\otimes Y:=(1\tens{}Y)\otimes M
\end{equation*}
whenever $X\in\mc{C}$, $Y\in\mc{D}$ and $M\in\mc{M}$. In the sequel we will not hesitate to work with either the left or right module category structures of a bimodule category without warning. 

Conversely, any pair of left $\mc{C}$-module and right $\mc{D}$-module category structures on $\mc{M}$ corresponds to a $\mc{C},\mc{D}$-bimodule category structure provided there is a family of natural isomorphisms $\gamma_{X, M, Y}:(X\otimes M)\otimes Y\simeq X\otimes(M\otimes Y)$ satisfying certain diagrams (see \cite[Proposition 2.12]{JG}). This bimodule category structure is strict provided both left and right module structures are strict and $\gamma$ is trivial. 
\end{defn}
\begin{defn}\label{moduleFunctorDef} Let $\mc{M}, \mc{N}$ be left $\mc{C}$-module categories and let  $F:\mc{M}\rightarrow\mc{N}$ be a functor. Then the pair $(F, f)$ is said to be a $\mc{C}$-\textit{module functor} if $f_{X, M}:F(X\otimes M)\rightarrow X\otimes F(M)$ are natural isomorphisms satisfying the usual coherence diagrams (again see \cite{O:MC}). We will often write just $F$ and leave $f$ implicit. A natural transformation $\tau:F\Rightarrow G$ for bimodule functors $(F, f), (G, g):\mc{M}\rightarrow\mc{N}$ is said to be a \textit{module natural transformation} whenever the diagram
\[
\xymatrix{
F(X\otimes M)\ar[rr]^{\tau_{X\otimes M}}\ar[d]_{f_{X, M}}&&G(X\otimes M)\ar[d]^{g_{X, M}}\\
X\otimes F(M)\ar[rr]_{id_X\otimes \tau_{M}}&&X\otimes G(M)
}
\]
commutes for all $X\in\mc{C}$ and $M\in\mc{M}$.
\end{defn}
\begin{notation} We will signify the 2-category of $(\mc{C}, \mc{D})$-bimodule categories having 1-cells bimodule functors, 2-cells bimodule natural transformations, by $(\mc{C}, \mc{D})$-Bimod. If $\mc{C}=\mc{D}$ we will write $\mc{C}$-Bimod.
\end{notation}
Unless otherwise stated we will require that all module categories be \textit{exact} in the sense of \cite{EO:FTC}. In particular we will assume that all module functors are exact. This is primarily a convenience for us since we will be working extensively with the relative tensor product (defined below).

With module categories, as with monoidal categories, it is convenient to require that the associativity and unit constraints are trivial, i.e. that the module structure is \textit{strict}. Fortunately this may always be assumed as the next theorem shows. Its proof is a straightforward extension of the proof of the MacLane strictness theorem for monoidal categories given in \cite{JS}. We state the result here without proof, all details of which may be found in the author's thesis \cite{JGthesis}.
\begin{thm}\label{MacLaneModule} Any module category is module equivalent to a strict module category.
\end{thm}
\subsection{2-categories}\label{2-catsubsect} In the sequel we will have occasion to make calculations in the 2-category setting. Recall that for 2-categories we must not only specify 0-cells (objects) and 1-cells (morphisms) as for usual categories, but also 2-cells (morphisms between morphisms) together with composition rules dictating how cells at different levels are to interact. We adopt the following notation: composition between 1-cells $f, g$ will be denoted $f\circ g$ or by juxtaposition when such compositions are defined. The \textit{vertical composition} of 2-cells $\sigma, \tau:f\rightarrow g$ will be denoted $\sigma\circ\tau:f\rightarrow g$. For 2-cells $\alpha:f\rightarrow g$, $\beta:f'\rightarrow g'$ where each 1-cell $f', g'$ is composable with each of $f, g$ denote their \textit{horizontal composition} by $\beta*\alpha:f'f\rightarrow g'g$. Composition between cells of different degrees will also be denote by $*$. For example the 2-cell that results by composing 2-cell $\alpha:f\rightarrow g:X\rightarrow Y$ with one cell $h:Y\rightarrow Z$ will be denoted $h*\alpha:hf\rightarrow hg$. 

These various compositions are required to satisfy certain rules of engagement. As an example we must have $(kh)*\alpha=k*(h*\alpha):khf\rightarrow khg$ for $\alpha, h$ above and 1-cell $k:Z\rightarrow W$. Also vertical and horizontal compositions must satisfy
\begin{equation}\label{2cellcompcompat}
(g*\beta)\circ(\alpha*f')=(\alpha*g')\circ(f*\beta)
\end{equation}
In general composition of 2-cells in a 2-category may take the form of polytopes with vertices 0-cells, edges 1-cells and faces labelled by 2-cells. For example the pair of planar diagrams below form a polytope (in this case a cube) when pasted together along the bold perimeters where they agree. 
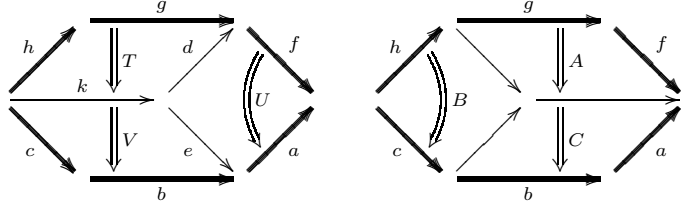
\[
\xymatrix{
& \arb{rr}^<<<<{}="1"^g & & \arb{dr}_<<<<{}="5"^f &\\
\arb{ur}^h \ar[rr]^>>>>>>{}="2"_>>>>>>{}="3"^k \arb{dr}_c & &\ar[ur]^d \ar[dr]_e & &\\
& \arb{rr}^<<<<{}="4"_b & &\arb{ur}^<<<<{}="6"_a &\\
\ar@{=>}"1";"2"^{T} 
\ar@{=>}"3";"4"^{V} 
\ar@/_.6pc/@{=>}"5";"6"^{U}
}\;\quad 
\xymatrix{
&\arb{rr}^>>>>>>{}="1"^g\ar[dr] & & \arb{dr}^f &\\
\arb{ur}_>>>>{}="5"^h \arb{dr}^>>>>{}="6"_c & &  \ar[rr]^<<<<{}="2"_<<<<{}="3" & &\\
& \arb{rr}^>>>>>>{}="4"_b  \ar[ur]& &\arb{ur}_a &\\
\ar@{=>}"1";"2"^{A} 
\ar@{=>}"3";"4"^{C} 
\ar@/^.6pc/@{=>}"5";"6"^{B}
}
\vspace{-2.5em}\]
Equality of the 2-compositions appearing in the planar diagrams is equivalent to saying that the associated polytope commutes. See \cite{BenBicat}, \cite{HigherOpsCats} or \cite{KV} for further details. 

In this paper the only 2-categories which will appear are those with 0-cells categories of some specified type, 1-cells functors and 2-cells natural transformations. In this case $\circ$-compositions of unmixed type are the usual compositions of functors and transformations. Let $\alpha:F\rightarrow G:\mc{B}\rightarrow\mc{C}$ be a natural transformation between functors $F, G$. Then for functors $E:\mc{A}\rightarrow\mc{B}$, $H:\mc{C}\rightarrow\mc{D}$ the compositions $\alpha*E:FE\rightarrow GE$ and $H*\alpha:HF\rightarrow HG$ are 2-cells having components $(\alpha*E)_X:=\alpha_{E(X)}$ and $(H*\alpha)_Y:=H(\alpha_Y)$ for $X\in\mc{A}$ and $Y\in\mc{B}$, respectively. 
\subsection{Balancing} In this section we prepare the ground for the definition of the relative tensor product by providing basic but necessary definitions. For further  discussion see \cite{ENO:homotop} or \cite{JG}.
\begin{defn}\label{DEF;bal} 
%
Let $\mc{C}$ be a monoidal category. Suppose $(\mc{M}, \mu)$, $(\mc{N}, \eta)$ are left, right $\mc{C}$-module categories having module structures $\mu, \eta$ respectively. Then the pair $(F, b)$ is called a \textit{$\mc{C}$-balanced functor} whenever $F:\mc{M}\boxtimes \mc{N} \rightarrow \mc{A}$ is a functor and $b_{M, X, Y}:F((M\otimes X)\boxtimes N)\simeq F(M\boxtimes (X\otimes N))$ are natural isomorphisms satisfying the pentagon 
\[\!\!\!\!\!\!\!\!\!\!\!\!\!\!
\xymatrix{
F((M\otimes(X\otimes Y))\boxtimes N) \ar[rr]^{b_{M, X\otimes Y, N}} \ar[d]^{\mu_{M, X, Y}}& & F(M\boxtimes((X\otimes Y)\otimes N)) \ar[d]^{\eta_{X, Y, N}}\\
F(((M\otimes X)\otimes Y)\boxtimes N) \ar[dr]_{b_{M\otimes X, Y, N}} & & F(M\boxtimes(X\otimes (Y\otimes N))) \\
& F((M\otimes X)\boxtimes (Y\otimes N)) \ar[ur]_{b_{M, X, Y\otimes N}} }
\]
for any $X$, $Y\in\mc{C}$ and $M\in\mc{M}$. Here the natural isomorphism $b$ will be referred to as the balancing morphism for the functor $F$. 

Put differently, $F$ is $\mc{C}$-balanced whenever $b:F(\otimes\tens{} 1)\stackrel{\sim}{\rightarrow} F(1\tens{}\otimes):\mc{M}\tens{}\mc{C}\tens{}\mc{N}\rightarrow\mc{A}$ is a natural isomorpism making the polytope commute on the level of 2-cells:
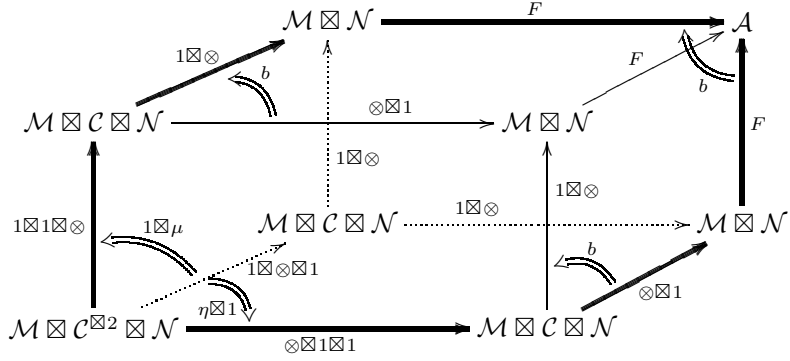
\[
\xymatrix{
& \mc{M}\tens{}\mc{N}\arb{rr}_>>>>>>{}="11"^F& &\mc{A}\\
\mc{M}\tens{}\mc{C}\tens{}\mc{N}\arb{ur}_<<<<<<<<<<<{\,\,\,}="1"^{1\tens{}\otimes} \ar[rr]^>>>>>>>>>>>>>>{\otimes\tens{}1}^<<<<<<<<<<<<<<{}="2"_>>>>>>{}="9" &&\mc{M}\tens{}\mc{N}\ar[ur]^F &\\
&\mc{M}\tens{}\mc{C}\tens{}\mc{N}\ar@{.>}[uu]|{\,\,\,}_<<<<<<<{1\tens{}\otimes}\ar@{.>}[rr]|{\,\,}^<<<<<<<<<<{1\tens{}\otimes} & &\mc{M}\tens{}\mc{N}\arb{uu}^>>>>>>{}="12" _F\\
\mc{M}\tens{}\mc{C}^{\tens{}2}\tens{}\mc{N}\ar@{.>}[ur]_<<<<<<<<<{}="3"_<<<<<<<<<{}="5"_>>>>{1\tens{}\otimes\tens{}1}\arb{uu}_<<<<<<<<<{}="4" ^{1\tens{}1\tens{}\otimes}\arb{rr}^<<<<<<<<<{}="6"_{\otimes\tens{}1\tens{}1}&&\mc{M}\tens{}\mc{C}\tens{}\mc{N}\ar[uu]_<<<<<<{}="7"^>>>>>>{}="10"_>>>>>>>{1\tens{}\otimes}\arb{ur}^<<<<<<{}="8"_{\otimes\tens{}1} &
\ar@/_.6pc/@{=>}"2";"1"_{b}
\ar@/_.6pc/@{=>}"3";"4"_{1\tens{}\mu}
\ar@/^.6pc/@{=>}"5";"6"_>>>{\eta\tens{}1}
\ar@/_.6pc/@{=>}"8";"7"_{b}
\ar@/^.6pc/@{=>}"12";"11"^>>>>>>{b}
}\]
The pentagon appearing at the beginning of this definition expresses commutativity of the polytope above but on the level of components (the associated 1-cells). 
For balanced functors $F, G:\mc{M}\boxtimes \mc{N} \rightarrow \mc{A}$ a natural transformation $\tau:F\rightarrow G$ is said to be balanced if it satisfies the diagram
\begin{equation*}
g_{M, X, N}\tau_{M\otimes X, N}=\tau_{M, X\otimes N}f_{M, X, N}
\end{equation*}
for all $M\in \mc{M}$, $N\in\mc{N}$ and $X$ in $\mc{C}$. Here $f, g$ are balancing morphisms for $F, G$. 
\end{defn}
A $\mc{C}$-balanced functor $(G, g):\mc{M}\rightarrow\mc{N}$ composes with any functor $F:\mc{N}\rightarrow\mc{P}$ yielding a $\mc{C}$-balanced functor $F\circ(G, g)=(F\circ G, F* g)$ ($*$ is as in \S\ref{2-catsubsect}). Balanced functors $\mc{M}\tens{}\mc{N}\rightarrow\mc{A}$ form a category having 1-cells given by balanced natural transformations. For further discussion see \S 3 of \cite{JG}.
\begin{remark}\label{C-CentralFunctor} It is possible to generalize Definition \ref{DEF;bal} slightly by defining balanced functors in the following way. 
\begin{defn}\label{BalanceGeneralDefinition} Let $(\mc{M}, \mu)$ be a $\mc{C}$-bimodule category and let $\mc{A}$ be any abelian category. Then the pair $(F, b)$ is called a $\mc{C}$-\textit{balanced functor} (or is said to be $\mc{C}$-balanced \textit{with respect to} $\mu$) whenever $F:\mc{M}\rightarrow\mc{A}$ is a functor and $b:F\otimes\rightarrow F\otimes(12):\mc{M}\tens{}\mc{C}\rightarrow\mc{A}$ is a natural isomorphism such that the following polytope commutes on the level of 2-cells:
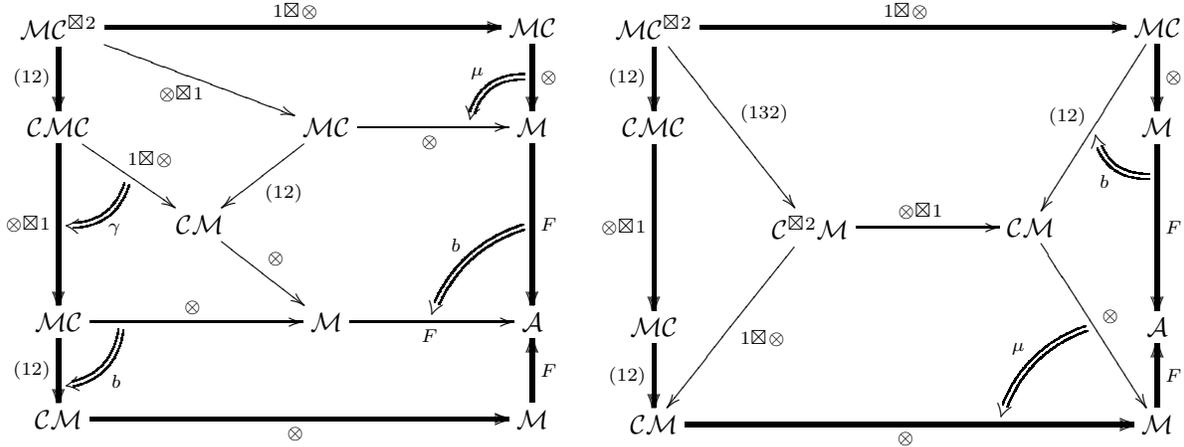
\[\!\!\!\!\!\!\!\!\!\!\!\!\!\!\!\!\!\!\!\!\!\!\!\!\!\!\!\!
\xymatrix{
\mc{M}\mc{C}^{\tens{}2}\arb{rrrr}^{1\tens{}\otimes}\ar[drr]_{\otimes\tens{}1}\arb{d}_{(12)}&&&&\mc{M}\mc{C}\arb{d}^\otimes_{}="7"\\
\mc{C}\mc{M}\mc{C}\arb{dd}_{\otimes\tens{}1}^{}="2"\ar[dr]^{1\tens{}\otimes}_{}="1"&&\mc{M}\mc{C}\ar[rr]_\otimes^>>>>>>{}="8"\ar[dl]^{(12)}&&\mc{M}\arb{dd}^F_{}="6"\\
&\mc{C}\mc{M}\ar[dr]^\otimes&&&\\
\mc{M}\mc{C}\arb{d}_{(12)}^>>>{}="4"\ar[rr]^\otimes_<<<<<{}="3"&&\mc{M}\ar[rr]_F^{}="5"&&\mc{A}\\
\mc{C}\mc{M}\arb{rrrr}_\otimes&&&&\mc{M}\arb{u}_F&
\ar@/^.6pc/@{=>}"1";"2"^{\gamma}
\ar@/^.6pc/@{=>}"3";"4"^{b}
\ar@/_.6pc/@{=>}"6";"5"_{b}
\ar@/_.6pc/@{=>}"7";"8"_{\mu}
}\!\!\!\!\!\!\!\!\!\!\!\!\!\;
\xymatrix{
\mc{M}\mc{C}^{\tens{}2}\arb{rrrr}^{1\tens{}\otimes}\ar[ddr]^{(132)}\arb{d}_{(12)}&&&&\mc{M}\mc{C}\arb{d}^\otimes\ar[ddl]_{(12)}^{}="10"&\\
\mc{C}\mc{M}\mc{C}\arb{dd}_{\otimes\tens{}1}&&&&\mc{M}\arb{dd}^F_<<<<<{}="9"\\
&\mc{C}^{\tens{}2}\mc{M}\ar[ddl]^{1\tens{}\otimes}\ar[rr]^{\otimes\tens{}1}&&\mc{C}\mc{M}\ar[ddr]^\otimes_{}="11"&\\
\mc{M}\mc{C}\arb{d}_{(12)}&&&&\mc{A}\\
\mc{C}\mc{M}\arb{rrrr}_\otimes^>>>>>>>>>>>>>>>>>>{}="12"&&&&\mc{M}\arb{u}_F
\ar@/^.6pc/@{=>}"9";"10"^{b}
\ar@/_.6pc/@{=>}"11";"12"_{\mu}
}
\]
Here we have abbreviated Deligne product as juxtaposition, $\gamma$ is the bimodule consistency constraint provided as part of the bimodule structure in $\mc{M}$, $\mu$ refers to either the left or the right $\mc{C}$-module category structure coming from the bimodule constraint, and we have indicated the various twists by their acting permutations. On the level of components commutativity of the polytope becomes the hexagon
\[\xymatrix{
F((X\otimes Y)\otimes M)\ar[d]_{\mu_{X, Y, M}}\ar[rr]^{b_{X\otimes Y, M}}&&F(M\otimes(X\otimes Y))\ar[d]^{\mu_{M, X, Y}}\\
F(X\otimes(Y\otimes M))\ar[d]_{b_{X, Y\otimes M}}&&F((M\otimes X)\otimes Y)\\
F((Y\otimes M)\otimes X)\ar[rr]_{\gamma_{Y, M, X}}&&F(Y\otimes (M\otimes X))\ar[u]_{b_{Y, M\otimes X}}
}\]
\end{defn}
%

To see that this new definition generalizes Definition \ref{DEF;bal} let $\mc{M}_i$, $1\leq i\leq n$, be $\mc{C}$-bimodule categories and define a left $\mc{C}$-module structure \textit{in the $i^{\textrm{th}}$ place} on $\mc{M}:=\mc{M}_1\tens{}\cdots\tens{}\mc{M}_n$ by $\otimes_i:=\otimes\tau_i$ where $\otimes$ is the left $\mc{C}$-module structure in $\mc{M}_i$ and $\tau_i$ is the twist functor operating via the permutation $(i, i-1,\dots,2,1)$:
\begin{equation*}
\tau_i:\mc{C}\tens{}\mc{M}\rightarrow \mc{M}_1\tens{}\cdots\tens{}\mc{M}_{i-1}\tens{}\mc{C}\tens{}\mc{M}_i\tens{}\cdots\tens{}\mc{M}_n.
\end{equation*}
Thus, for $X\in\mc{C}$ we have
\begin{equation*}
X\otimes_i(M_1\tens{}\cdots\tens{}M_n):=M_1\tens{}\cdots M_{i-1}\tens{}(X\otimes M_i)\tens{}\cdots\tens{} M_n.
\end{equation*}
Similarly define right $\mc{C}$-module category structure $\otimes^j$ \textit{in the $j^{\textrm{th}}$ place} using the right $\mc{C}$-module category structure in $\mc{M}_j$ and a general twist $\tau^j$ acting from $\mc{M}\tens{}\mc{C}$: 
\begin{equation*}
(M_1\tens{}\cdots\tens{}M_n)\otimes^jX:=M_1\tens{}\cdots M_{j-1}\tens{}(M_j\otimes X)\tens{}\cdots\tens{} M_n.
\end{equation*}
The module category structures $\otimes^j$ and $\otimes _{i}$ give $\mc{M}$ the structure of a $\mc{C}$-bimodule category. Denote by $\mc{M}^{i,j}$ the category $\mc{M}$ with this bimodule structure. The bimodule consistency constraint ($\gamma$ in Definition \ref{BimoduleCatDef}) in $\mc{M}^{i,i}$ is precisely that in $\mc{M}_i$, and is trivial in the case $i\neq j$. Define a functor $F:\mc{M}\rightarrow\mc{A}$ to be $(i,j)$\textit{-balanced with respect to $\mc{C}$} if it is balanced with respect to the bimodule structure $\mc{M}^{i,j}$ in the sense of Definition \ref{BalanceGeneralDefinition}. Then we may specialize and say that $F$ is balanced in position $j$ if it is $(j+1, j)$-balanced with respect to $\mc{C}$. Thus a functor from $\mc{M}\tens{}\mc{N}$ is $\mc{C}$-balanced as in Definition \ref{DEF;bal} if it is balanced in position 1. 

We can define what it means for a functor to be balanced in more than one position as follows. Suppose $|i - j|\geq1$ and take $\mc{M}$ as above. Then a functor $F:\mc{M}\rightarrow\mc{A}$ is $\mc{C}$-\textit{balanced in the $i^{\textrm{th}}$ and $j^{\textrm{th}}$ positions simultaneously} provided that it is balanced in position $i$ and balanced in position $j$, and the balancing isomorphisms $b^i:F(\otimes^i)\rightarrow F(\otimes_{i+1})$ and $b^j:F(\otimes^j)\rightarrow F(\otimes_{j+1})$ satisfy, in addition to the poytopes appearing in Definition \ref{BalanceGeneralDefinition} (one for each index), a balancing consistency condition given by a certain polytope. We won't need the general formulation\footnote{The general polytope is very similar to the one appearing below, but with a (yet to be defined) generalized balancing associator replacing $\gamma$ in various cells.} so we give this polytope only in the case $i=1$, $j=2$ below, abbreviating the Deligne product of abelian categories by juxtaposition.  
\[\!\!\!\!\!\!\!\!\!\!\!\xymatrix{
\mc{M}\mc{C}\mc{C}\arb{r}^{(123)}\arb{d}_{\otimes^1\tens{}1}&\mc{C}\mc{M}\mc{C}\arb{rr}^{\otimes_3}\ar[d]^{\otimes^3}&&\mc{M}\mc{C}\arb{dd}^{(12)}_{}="4"\ar[ddl]_{\otimes^1}^{}="3"\\
\mc{M}\mc{C}\ar[r]_{(12)}\arb{dd}_{\otimes^2}&\mc{C}\mc{M}\ar[dr]^{\otimes_3}&&\\
&&\mc{M}\ar[dl]^F_{}="2"&\mc{C}\mc{M}\arb{d}^{\otimes_2}\\
\mc{M}\arb{r}_F^{}="1"&\mc{A}&&\mc{M}\arb{ll}^F
\ar@/^1pc/@{=>}"1";"2"^{b^2}
\ar@/_.4pc/@{=>}"3";"4"_{b^1}
}\;\quad
\xymatrix{
\mc{M}\mc{C}\mc{C}\arb{r}^{(12)}\arb{d}_{\otimes^2}&\mc{C}\mc{M}\mc{C}\arb{rr}^{\otimes_2}_{}="10"\ar[d]_{\otimes^3}^{}="9"&&\mc{M}\mc{C}\arb{dd}^{(12)}_{}="8"\ar[ddl]_{\otimes^2}^{}="7"\\\
\mc{M}\mc{C}\ar[r]_{(12)}\arb{dd}_{\otimes^1}&\mc{C}\mc{M}\ar[dr]^{\otimes_2}&&\\
&&\mc{M}\ar[dl]_{}="6"^F&\mc{C}\mc{M}\arb{d}\\
\mc{M}\arb{r}_F^{}="5"&\mc{A}&&\mc{M}\arb{ll}^F
\ar@/^1pc/@{=>}"5";"6"^{b^1}
\ar@/_.4pc/@{=>}"7";"8"_{b^2}
\ar@/_.8pc/@{=>}"9";"10"_{\gamma}
}\]
%
In order to simplify the exposition, and because they commute trivially, we have not included the rectangular 2-cell with edges given by the pair of edge compositions $\mc{M}\mc{C}\mc{C}\rightarrow\mc{M}\mc{C}\rightarrow\mc{M}$, one from each of the perimeters of the above figures, nor have we included the octagonal 2-cell with edges formed by the two different edge compositions $\mc{M}\mc{C}\mc{C}\rightarrow\mc{C}\mc{M}\mc{C}\rightarrow\mc{M}\mc{C}\rightarrow\mc{C}\mc{M}\rightarrow\mc{M}$. One can check easily that each of the pairs of compositions agree. Written in linear equation form commutativity of the polytope above means
\begin{equation*}\label{MultiBalancedCondition1}
(b^1*\otimes_2)(F*\gamma)(b^1*\otimes^2) = (b^1*\otimes_3)(b^2*\otimes^1)
\end{equation*}
which becomes, on the level of components,
\begin{equation}\label{MultiBalancedCondition2}
b^2_{R, XS, Y, T}\circ F(1_R\tens{}\gamma_{Y, S, X}\tens{}1_T)\circ b^1_{R, X, SY, T}=b^1_{R, X, S, YT}\circ b^2_{RX, S, Y, T}
\end{equation}
for $X, Y\in\mc{C}$ and $R\tens{}S\tens{}T\in\mc{M}$. Equation (\ref{MultiBalancedCondition2}) appears as the diagram in Definition 3.4 in \cite{JG}. For $i = j=1$ and $n=2$ the associated polytope reduces to a diagram similar to that appearing in the proof of Proposition 4.11 in \textit{loc. cit}. We leave details to the interested and motivated reader. In a future article the phenomenon of balancing will be discussed in general at greater length. 
%
\end{remark}
\subsection{Relative tensor product} Just as the classical tensor product of modules is defined as a certain universal object for middle-balanced morphisms we may define the tensor product of module categories to be the abelian category universal with respect to balanced functors.
\begin{defn}\label{DEF;tensor} The \textit{tensor product} of right $\mc{C}$-module category $\mc{M}$ and left $\mc{C}$-module category $\mc{N}$ consists of an abelian category $\mc{M}\boxtimes_{\mc{C}}\mc{N}$ and a right exact $\mc{C}$-balanced functor $B_{\mc{M}, \mc{N}}: \mc{M}\boxtimes \mc{N}\rightarrow \mc{M}\boxtimes_{\mc{C}}\mc{N}$ universal for \textit{right exact} $\mc{C}$-balanced functors from $\mc{M}\boxtimes \mc{N}$. That is, for any balanced functor $F$ (squiggly arch in diagram below) there is a unique $\overline{F}$ with the property that $\overline{F}\Buni{M}{N}=F$. 
\end{defn}
\begin{remark}\label{universalNatTransRemark}
The category $\bitens{M}{C}{N}$ is unique up to a unique equivalence: if $U_{\mc{M}, \mc{N}}:\mc{M}\tens{}\mc{N}\rightarrow\mc{U}$ is another functor universal for balanced functors from $\bitens{M}{}{N}$ then there is a unique equivalence $A:\mc{U}\stackrel{\sim}{\rightarrow}\bitens{M}{C}{N}$ satisfying $AU_{\mc{M}, \mc{N}}=\Buni{M}{N}$. 
\end{remark}
\begin{defn}\label{DescendedBalancedNatTrans}Let $\tau:F\rightarrow G$ be a balanced natural transformation. Then there is a unique natural transformation $\overline{\tau}:\overline{F}\rightarrow\overline{G}$ having the property that $\overline{\tau}*\Buni{M}{N}=\tau$. This situation is illustrated in the following diagram of 2-cells where we understand each sub-triangle of 1-cells with common edge $\Buni{M}{N}$ (one taken with two squiggly arches, the other with all plain ones) to commute separately.
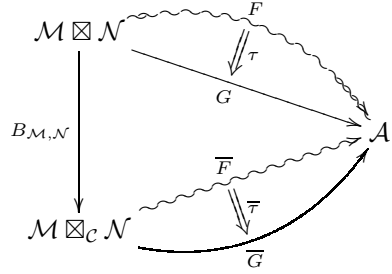
\[\xymatrix{
\mc{M}\tens{}\mc{N}\ar@/^2pc/@{~>}[drrr]^{F}_{}="1"\ar[drrr]_{G}^{}="2"\ar[dd]_{\Buni{M}{N}}&&&
\\
&&&\mc{A}\\
\mc{M}\tens{C}\mc{N}\ar@{~>}[urrr]^{\overline{F}}_{}="3"\ar@/_2pc/[urrr]_{\overline{G}}^{}="4"&&&
\\
\ar@{=>}"1";"2"^{\tau} \ar@{=>}"3";"4"^{\overline{\tau}} 
}
\vspace{-2.5em}\]
Using the universality of $\Buni{M}{N}$ and basic relations in the 2-category of functors and natural transformations one shows easily that, for composable balanced natural transformations $\tau$ and $\sigma$, $\overline{\tau\sigma}=\overline{\tau}\;\overline{\sigma}$.
\end{defn}
\begin{remark}\label{multi-relative tensor} Let $\mc{M}$ be as in Remark \ref{C-CentralFunctor}, with $\mc{C}$-bimodule structures at the $i^{\textrm{th}}$ and $j^{\textrm{th}}$ places simultaneously as described there. Then the simultaneous relative tensor product 
\begin{equation*}
B^{ij}:\mc{M}\rightarrow\mc{M}_1\tens{}\cdots\tens{}\mc{M}_i\tens{C}\mc{M}_{i+1}\tens{}\cdots\tens{}\mc{M}_j\tens{C}\mc{M}_{j+1}\tens{}\cdots\mc{M}_n
\end{equation*}
is the unique abelian category factoring through functors simultaneously $\mc{C}$-balanced in the $i^{\textrm{th}}$ and $j^{\textrm{th}}$ places. The universal functor $B^{ij}$ identifies with the Deligne product of functors $B^{ij}=id_1\tens{}\cdots\tens{}B_{i, i+1}\tens{}\cdots\tens{}B_{j, j+1}\tens{}\cdots\tens{}id_n$ where we have abbreviated $B_{s, t}$ as the universal $\mc{C}$-balanced functor from $\mc{M}_s\tens{}\mc{M}_t$. A natural transformation is balanced in the $i^{\textrm{th}}$ and $j^{\textrm{th}}$ places simultaneously if it is balanced in each separately (see Definition \ref{DEF;bal}).
\end{remark}
\subsection{Module centers.} It is possible to extend the notion of the center of a monoidal category to module categories in such a way that the new construction reduces to the monoidal center in the regular module category case. We review the construction here and refer to \cite{ENO:homotop}, \cite{JG}, \cite{GNN} where monoidal category centers appear in various contexts.
\begin{defn}\label{defCent} Let $\mc{M}$ be a strict $\mc{C}$-bimodule category. The \textit{relative center} $\cent{M}{C}$ of $\mc{M}$ over $\mc{C}$ consists of objects given by pairs $(M, \varphi_M)$ where $M\in\mc{M}$ and where, for any $X\in\mc{C}$,  $\varphi_{X, M}:X\otimes M\simeq M\otimes X$ are isomorphisms natural in $\mc{C}$ satisfying the diagram below for any $Y\in\mc{C}$.
\[
\xymatrix{
X\otimes M\otimes Y\ar[rr]^{\varphi_{X, M}\otimes id_Y}&&M\otimes X\otimes Y\\
&X\otimes Y\otimes M\ar[ul]^{id_X\otimes \varphi_{Y, M}}\ar[ur]_{\varphi_{X\otimes Y, M}}
}
\]
Call $\varphi$ a \textit{central structure} for $\mc{M}$ over $\mc{C}$. A morphism from $(M, \varphi_M)$ to $(N, \varphi_N)$ in $\cent{M}{C}$ is a morphism $t:M\rightarrow N$ in $\mc{M}$ satisfying $\varphi_{X, N}(id_X\otimes t) = (t\otimes id_X)\varphi_{X, M}$.
\end{defn}
\begin{note}\label{non-strictRelativeCenterNote} In the case that $\mc{M}$ is \textit{not} strict as a $\mc{C}$-bimodule category the diagram in Definition \ref{defCent} is no longer a triangle but a hexagon wherein appear left and right associativity (for the underlying left and right module category structures) as well as the bimodule consistency isomorphism $\gamma$ appearing in Definition \ref{BimoduleCatDef}. In case the left and right module structures are both strict but the bimodule consistency constraint isn't trivial the hexagonal diagram reduces to the rectangle
\[
\xymatrix{
X\otimes Y\otimes M\ar[rr]^{\varphi_{X\otimes Y, M}}\ar[d]_{id_X\otimes\varphi_{Y, M}}&&M\otimes X\otimes Y\\
X\otimes(M\otimes Y)\ar[rr]_{\gamma_{X, M, Y}}&&(X\otimes M)\otimes Y\ar[u]_{\varphi_{X, M}\otimes id_Y}
}
\]
\end{note}
It is known that $\cent{M}{C}$ has canonical structure of a $\mc{Z}(\mc{C})$-module category and that, as such, it is equivalent to a certain category of right exact $\mc{C}$-bimodule functors. We reproduce the statement here in full. 
\begin{thm}[\cite{JG}]\label{centralproplem} For $\mc{C}$-bimodule category $\mc{M}$ there is a canonical $\mc{Z}(\mc{C})$-bimodule equivalence $Fun_{\mc{C}\boxtimes\mc{C}^{rev}}(\mc{C}, \mc{M})\simeq\cent{M}{C}$.
\end{thm}
\begin{proof}
Any $\mc{C}$-bimodule functor $F:\mc{C}\rightarrow \mc{M}$ is completely specified by its action on the unit object $1$ of $\mc{C}$. If $F(1)=M$ write $F:=F_M$. Then the equivalence $Fun_{\mc{C}\boxtimes\mc{C}^{rev}}(\mc{C}, \mc{M})\stackrel{\sim}{\rightarrow}\cent{M}{C}$ is defined by sending 
\begin{equation}\label{FunCentEquiv}
(F_M, f)\mapsto(M, \hat{f}_M),\quad\hat{f}_{X, M}:=f_{1\tens{}X, 1}\circ f^{-1}_{X\tens{}1, 1}:X\otimes M\rightarrow M\otimes X.
\end{equation}
\end{proof}
\begin{remark}\label{CenterExactModuleRemarkI} Theorem \ref{centralproplem}, together with \cite[Lemma 3.30]{EO:FTC}, implies that $\cent{M}{C}$ is exact as a module category over $\cent{C}{}$.
\end{remark}
\begin{remark}\label{ForgetFulCenter}
The obvious canonical forgetful functor $\Forg:\mc{Z}_\mc{C}(\mc{M})\rightarrow\mc{M}$ is defined by sending $(M, \varphi_M)\mapsto M$. Composing this with the equivalence above gives a forgetful functor from the category of $\mc{C}$-bimodule functors $\mc{C}\rightarrow\mc{M}$ defined by $F\mapsto F(1)$.
\end{remark}
In \cite{ENO:homotop} it is shown that, for $\mc{C}$ a tensor category and $\mc{M}$ a semisimple $\mc{C}$-bimodule category, the relative center $\mc{Z}_\mc{C}(\mc{M})$ is equivalent to the relative tensor product. Theorem \ref{centertensorthm} below extends this result to the exact case. 
\subsection{Pre-metric groups}\label{PreMetGrpSec} Everything in this subsection may be found in \cite{braidedI}. Additionally we refer the reader to \cite{K:QG} and \cite{BK} for definitions and other information relating to braided fusion categories.

Let $G$, $B$ be abelian groups. A \textit{quadratic form} on $G$ having values in $B$ is a map $q:G\rightarrow B$ satisfying $q(g^{-1})=q(g)$ such that the symmetric function $b(g, h):=\frac{q(gh)}{q(g)q(h)}$ is bimultiplicative. We call $b:G\times G\rightarrow B$ the bimultiplicative form associated to $q$. If $B=k^{\times}$ then we call $b$ the bicharacter associated to $q$.
\begin{defn}
A \textit{pre-metric group} is a pair $(G, q)$ where $G$ is a finite abelian group and $q:G\rightarrow k^\times$ is a quadratic form. A morphism of pre-metric groups $(G_1, q_1)\rightarrow (G_2, q_2)$ is a group homomorphism $f:G_1\rightarrow G_2$ satisfying the triangle $q_2\circ f=q_1$.  
\end{defn}
The set of isomorphism classes of the objects of any pointed braided fusion category $\mc{C}$ forms a group $G$. For $g\in G$ denote by $q(g)\in k^\times$ the braiding $c_{X, X}\in\Aut(X\otimes X)$ where $X$ is in the isomorphism class $g$. Then the association $g\mapsto q(g)$ constitutes a quadratic form $G\rightarrow k^\times$. In this way $\mc{C}$ determines the pre-metric group $(G, q)$.

Conversely every pre-metric group $(G, q)$ determines a pointed braided fusion category $\mc{C}(G, q)$.  As a fusion category $\mc{C}(G, q)$ is $Vec_G$, the category of finite-dimensional $G$-graded vector spaces. For a homogeneous object $X$ of degree $g$ define the twist $\theta_X=q(g)$. Then the braiding $c_{X, Y}:X\otimes Y\rightarrow Y\otimes X$ satisfies the equation
\begin{equation}\label{ccdef}
c_{X, Y}c_{Y, X}=b(g, h)id_{Y\otimes X}
\end{equation}
where $b$ is the bicharacter determined by $q$. In the special case that $q$ comes from a bicharacter $\beta:G\times G\rightarrow k^\times$ via $q(x)=\beta(x,x)$ the associated braiding is $c_{X, Y}=\beta(g, h)\tau$ (here $\tau$ is the linear twist). These two constructions define reciprocal equivalences between the category of pre-metric groups and the (truncated 2-) category of pointed braided fusion categories. We refer the reader to the original paper of Joyal and Street (\cite{JS}) for further details. 
\subsection{Dominant functors.} Let $F:\mc{A}\rightarrow\mc{B}$ be an additive functor between abelian categories and define its image $Im(F)$ to be the full subcategory of $\mc{B}$ having objects given by all subquotients of objects of the form $F(X)$ for any $X\in\mc{A}$.
\begin{defn}\label{DominantFunctorDef} The functor $F$ is said to be \textit{dominant} if $Im(F)=\mc{B}$.
\end{defn}
It is an easy exercise to show that $Im(F)$ is itself an abelian category. Furthermore, if $\mc{A}, \mc{B}$ are tensor categories and $F$ a tensor functor then $Im(F)$ is a tensor subcategory of $\mc{B}$. Indeed, if $A_1, A_2$ are quotients of some subobjects $Z_1, Z_2$ of $F(X_1), F(X_2)$ resp., then exactness of the tensor structure in $\mc{B}$ implies that $A_1\otimes A_2$ is a quotient of $Z_1\otimes Z_2$, a subobject of $F(X_1)\otimes F(X_2)\simeq F(X_1\otimes X_2)$. Hence $A_1\otimes A_2$ can be identified with a subquotient of $F(X_1\otimes X_2)$ and therefore with an object of $Im(F)$. The unit object $1$ is contained in $Im(F)$ because it is a subobject of $F(1)$. All structural constraints come from those in $\mc{B}$. It is also evident that if $\mc{A}, \mc{B}$ are semisimple then dominance of $F$ means that any object of $\mc{B}$ is actually a subobject of $F(X)$ for some $X\in\mc{A}$. 

We will require a special case of the following lemma.
\begin{lem}\label{dominance of B Lemma} Let $\mc{C}$ be a tensor category, and let $\mc{M}, \mc{N}$ be right, left $\mc{C}$-module categories, respectively. Then the universal balanced functor $B_{\mc{M}, \mc{N}}$ appearing in Definition \ref{DEF;tensor} is dominant.
\end{lem}
\begin{proof}
Let $F$ be any balanced functor from $\mc{M}\tens{}\mc{N}$, and let $\overline{F}$ be the unique functor from $\mc{M}\tens{C}\mc{N}$ with $\overline{F}\Buni{M}{N}=F$. For the inclusion $i:Im(\Buni{M}{N})\hookrightarrow\mc{M}\tens{C}\mc{N}$ define $F':=\overline{F}i$. Then evidently $F'\Buni{M}{N}=F$, hence $F$ factors through $Im(\Buni{M}{N})$ uniquely, and as a consequence of the universality of the relative tensor product $\mc{M}\tens{C}\mc{N}=Im(\Buni{M}{N})$.
\end{proof}
\section{The relative center}\label{sectionRelCentMonCat}
In what follows all module categories are assumed to be exact, and we will not assume in general that tensor categories are semisimple. 

We begin this section by proving a convenient theorem illustrating a relationship between the relative center and the relative tensor product. The proof below essentially occurs in \cite[Section 5.8]{ENO:homotop} and we provide it here for completeness. Let $\mc{C}$ be a tensor category. 
\begin{thm}\label{centertensorthm}
Let $\mc{M}$ be a strict right $\mc{C}$-bimodule category and $\mc{N}$ a strict left $\mc{C}$-module category. Then there is a canonical equivalence of abelian categories $\cent{\bitens{M}{}{N}}{C}\simeq\bitens{M}{C}{N}$ such that the universal balanced functor $\Buni{M}{N}$ identifies with the canonical functor $\bitens{M}{}{N}\rightarrow \mc{Z}_{\mc{C}}(\bitens{M}{}{N})$. 
\end{thm}
\begin{proof} We begin by defining a certain functor.  
Let $\mc{A}$ be any exact $\mc{C}$-bimodule category. Consider the forgetful functor $Forg:\mc{Z}_\mc{C}(\mc{A})\rightarrow \mc{A}$. Since it has the structure of a module functor between exact module categories its left adjoint exists (thanks to \cite[Lemma 3.21]{EO:FTC}). Denote its left adjoint by $\mc{Z}_{\mc{A}, \mc{C}}:\mc{A}\rightarrow\cent{A}{C}$. Then $\mc{Z}_{\mc{A}, \mc{C}}$ is dominant and $\mc{C}$-balanced, and universal in the sense that if $Q:\mc{A}\rightarrow\mc{S}$ is any other $\mc{C}$-balanced functor there is a unique $\overline{Q}:\cent{A}{C}\rightarrow\mc{S}$ with $\overline{Q}\mc{Z}_{\mc{A}, \mc{C}}=Q$. Indeed, let $Q'$ be right adjoint to $Q$. Then for any $S\in\mc{S}$ the object $Q'(S)$ has the structure of an object of $\mc{Z}_{\mc{C}}(\mc{A})$. Furthermore, the functor $Q'':\mc{S}\rightarrow\mc{Z}_\mc{C}(\mc{A})$ sending $S\mapsto Q'(S)$ satisfies $Forg\circ Q'' = Q'$. Taking the left adjoint then gives $\overline{Q}\mc{Z}_{\mc{A}, \mc{C}}=Q$ where $\overline{Q}$ is the left adjoint of $Q''$. This is illustrated diagrammatically below.
\[
\xymatrix{
\mc{S}\ar[r]^{Q''}\ar[dr]_{Q'}&\mc{Z}_\mc{C}(\mc{A})\ar[d]^{\textrm{Forg}}="1"&&\mc{S}\ar@{}[d]_{}="2"&\mc{Z}_\mc{C}(\mc{A})\ar[l]_{\overline{Q}}\\
&\mc{A}&&&\mc{A}\ar[ul]^{Q}\ar[u]_{\mc{Z}_{\mc{A}, \mc{C}}}
\ar@{}"1";"2"^{\stackrel{\textrm{ad}}{\rightarrow}}
}\]

Let $\mc{M}$ be a right $\mc{C}$-module category, and $\mc{N}$ a left $\mc{C}$-module category. Then the Deligne product $\bitens{M}{}{N}$ has the structure of a $\mc{C}$-bimodule category via
\begin{equation}\label{canonicalonesideds}
X\otimes(M\tens{}N):=M\tens{}(X\otimes N),\qquad (M\tens{}N)\otimes X:=(M\otimes X)\tens{}N
\end{equation}
for any $X\in\mc{C}$, $M\in\mc{M}$ and $N\in\mc{N}$. From the beginning of this proof we have a $\mc{C}$-balanced functor $\mc{Z}_{\bitens{M}{}{N}, \mc{C}}:\bitens{M}{}{N}\rightarrow \cent{\bitens{M}{}{N}}{\mc{C}}$ universal for $\mc{C}$-balanced functors from $\bitens{M}{}{N}$. Using the stated bimodule structure in $\bitens{M}{}{N}$ this means this functor is universal for balanced functors from $\bitens{M}{}{N}$ (see Remark \ref{C-CentralFunctor}). 
Hence $\cent{\bitens{M}{}{N}}{\mc{C}}$ possesses the defining universal property unique to $\bitens{M}{C}{N}$ and must therefore be equivalent to it. Furthermore, the functor $\mc{Z}_{\bitens{M}{}{N}, \mc{C}}$ identifies with $B_{\mc{M}, \mc{N}}$ up to unique equivalence (see \cite[Remark 3.6]{JG}). 
\end{proof}
\begin{remark}
In the case that $\mc{A}=\mc{C}$ is fusion the functor $\mc{Z}_{\mc{C}, \mc{C}}:\mc{C}\rightarrow\mc{Z}_{\mc{C}}(\mc{C})=\mc{Z}(\mc{C})$ is the induction functor discussed in \cite[Section 5.8]{ENO:ofc}.
\end{remark}
%
%
Let $\mc{D}$ be a symmetric braided tensor category. Recall that a tensor category $\mc{C}$ is said to be tensor \textit{over} $\mc{D}$ if there is a braided inclusion $\mc{D}\hookrightarrow \mc{Z}(\mc{C})$ such that the composition $\mc{D}\hookrightarrow \mc{Z}(\mc{C})\rightarrow \mc{C}$ is faithful (Definition \ref{overnessdef}). In such a case we consider $\mc{D}$ to be a braided tensor subcategory of $\mc{C}$ by identifying it with its image therein. This inclusion induces in $\mc{C}$ the structure of a $\mc{D}$-bimodule category as described in the introduction. 
\begin{proof}[Proof of Theorem \ref{tensRelCentThm}]
The canonical tensor structure in $\cent{C}{D}$ is inherited directly from that in $\mc{C}$ by the formula
\begin{equation}\label{relcentTens}
(X,  c_X)\otimes (Y,  c_Y)=(X\otimes Y,  c_{X\otimes Y})
\end{equation}
where as usual $ c_{Z, X\otimes Y}:=(id_X\otimes  c_{Z, Y})( c_{Z, X}\otimes id_Y)$ for any $Z\in\mc{D}$. Observing that the category of $\mc{D}$-bimodule functors $Fun_{\mc{D}\tens{}\mc{D}}(\mc{D}, \mc{C})$ has tensor structure $F_X\otimes F_Y:=F_{X\otimes Y}$ the functor $\mc{Z}_{\mc{C}, \mc{D}}$ from the proof of Theorem \ref{centertensorthm} is evidently strict monoidal, being as it is the adjoint of a (strong) monoidal functor. 

To see that the relative center $\cent{C}{D}$ is tensor over $\mc{D}$ 
we must describe a braided inclusion $\mc{D}\hookrightarrow\mc{Z}(\cent{C}{D})$. For any $D\in\mc{D}$ and $(X,  c_X)\in\cent{C}{D}$ define the natural isomorphism
\begin{equation}\label{centerRelCentStruct}
\psi_{D, (X,  c_X)}:= c_{D, X}:D\otimes X\rightarrow X\otimes D.
\end{equation}
To see that $c_{D, X}$ is a morphism in $\cent{C}{D}$ for any $(X,  c_X)\in \cent{C}{D}$ observe commutativity of the diagram below. Let $D, E\in\mc{D}$.
\[\xymatrix{E\otimes D\otimes X\ar[drr]^{c_{E, D}\otimes id_X}\ar[dr]_{c_{E\otimes D, X}}\ar[dd]_{id_E\otimes c_{D, X}}\ar[rrr]^{c_{E, D\otimes X}}&&&D\otimes X\otimes E\ar[dd]^{c_{D, X}\otimes id_E}\\&X\otimes E\otimes D\ar[drr]_{id_X\otimes c_{E, D}}&D\otimes E\otimes X\ar[dr]^{c_{D\otimes E, X}}\ar[ur]^{id_D\otimes c_{E, X}}&\\E\otimes X\otimes D\ar[rrr]_{c_{E, X\otimes D}}\ar[ur]_{c_{E, X}\otimes id_X}&&&X\otimes D\otimes E
}\]
The scalene triangles are the definition of $c_{E, X\otimes D}$ and $c_{E, D\otimes X}$, the isosceles triangles are the diagrams defining $c$ as central structure, and the parallelogram is the square of naturality for $c$. 

Now define $\Psi:\mc{D}\rightarrow\mc{Z}(\cent{C}{D})$ by $D\mapsto(D, \psi_D)$ on objects and $\Psi=id$ on morphisms. Evidently $\Psi$ is an inclusion, and it therefore remains to show that it is braided. But this is trivial; $\mc{D}\hookrightarrow \mc{Z}(\mc{C})$ is braided, $\Psi(d_{E, D})=c_{E, D}$ and thus $\Psi$ is strictly braided

We now show that $\cent{C}{D}$ inherits rigidity from that in $\mc{C}$. For $X\in\mc{C}$ define $F_X^*:=F_{X^*}$ (here $F_A$ is the functor $X\mapsto X\otimes A$). The evaluation and coevaluation maps come from those in $\mc{C}$:
\begin{equation*}
\ev_{F_X}:=F_{\ev_X}:F^*_X\otimes F_X\rightarrow F_1,\quad \coev_{F_X}:=F_{\coev_X}:F_1\rightarrow F_X\otimes F^*_X.
\end{equation*}
Here we are writing $F_f:F_A\rightarrow F_B$ for the natural transformation coming from the morphism $f:A\rightarrow B$. The necessary compositions obtain for the pair $\ev_{F_X}, \coev_{F_X}$ since for composable morphisms $f, g$ one has $F_fF_g=F_{fg}$. Explicitly, for any $Z$, 
\begin{eqnarray*}
(1_{F_X}\otimes\ev_{F_X})_Z\circ(\coev_{F_X}\otimes 1_{F_X})_Z&=&(F_{1_X\otimes\ev_X})_Z\circ (F_{\coev_X\otimes 1_X})_Z\\
&=&(F_{(1_X\otimes\ev_X)(\coev_X\otimes 1_X)})_Z=(F_{1_X})_Z
\end{eqnarray*}
and since $(F_{1_X})_Z=(1_{F_X})_Z$ the initial composition is identity. The other identity required of rigidity is proved similarly. 

Finally suppose $\mc{C}$ is fusion. Using Theorem \ref{centralproplem} in conjunction with \cite[Theorem 2.16]{ENO:ofc} the category $\cent{C}{D}$ is semisimple. This completes the proof of  the theorem.
\end{proof} 
\begin{remark} With respect to the tensor structure inherited by the relative center the functor $\mc{Z}_{\mc{C}, \mc{D}}:\mc{C}\rightarrow\cent{C}{D}$ from the proof of Theorem \ref{centertensorthm} is strict monoidal. 
If $\mc{C}$ is fusion and $\mc{D}$ is semisimple then the relative center $\cent{C}{D}$ is also fusion.
\end{remark}
%
%
\begin{proof}[Proof of Corollary \ref{monoidalityOfTensor}] Let $\sigma_i:\mc{D}\hookrightarrow \mc{Z}(\mc{C}_i)$ be the braided inclusions putting $\mc{C}_i$ tensor over $\mc{D}$ for $i=1, 2$. These combine to form a braided inclusion $\sigma:=\sigma_1\tens{}\sigma_2:\mc{D}\tens{}\mc{D}\hookrightarrow \mc{Z}(\mc{C}_1)\tens{}\mc{Z}(\mc{C}_2)$. Using $\sigma$ we define the braided inclusions $\psi_i$ by precomposing $\sigma$ with the braided inclusion $X\mapsto X\tens{}1$, $X\mapsto 1\tens{}X$: 
\begin{equation*}
\psi_i:\mc{D}\hookrightarrow\mc{D}\tens{}\mc{D}\stackrel{\sigma}{\hookrightarrow}\mc{Z}(\mc{C}_1)\tens{}\mc{Z}(\mc{C}_2)\hookrightarrow\mc{Z}(\mc{C}_1\tens{}\mc{C}_2). 
\end{equation*}
Thus $\psi_1:D\mapsto \sigma(D\tens{}1)$, $\psi_2:D\mapsto\sigma(1\tens{}D)$ and $\mc{C}_1\tens{}\mc{C}_2$ is tensor over $\mc{D}$. Applying the equivalence $\mc{Z}_\mc{D}(\mc{C}_1\tens{}\mc{C}_2)\simeq \mc{C}_1\tens{D}\mc{C}_2$ from Theorem \ref{centertensorthm} and Theorem \ref{tensRelCentThm} shows that the relative tensor product has the structure of a tensor category.

In order for this tensor structure to be canonical we verify that the $\psi_i$ are braided equivalent. Abbreviate $D_1:=D\tens{}1$ and $D_2:=1\tens{}D$ for ay $D\in\mc{D}$ and write $\psi_1(D):=(D_1, c'_{D_1})$ and $\psi_2(D):=(D_2, c'_{D_2})$ where, as in Equation (\ref{TwistBraid}) below, $c'_{X\tens{}Y, U\tens{}V}:=c^1_{X, U}\tens{}c^2_{Y, V}$. Define the natural isomorphism $\tau:\psi_1\stackrel{\sim}{\rightarrow}\psi_2$ with components $\tau_D:\sigma_1(D)\tens{}1\simeq1\tens{}\sigma_2(D)$ by the image of the permutation action $(12)$ (just switch the order of factors).  It's trivial that $\tau_D$ is a morphism in $\mc{Z}(\mc{C}_1\tens{}\mc{C}_2)$. The following diagram shows that $\tau$ is a braided equivalence.
\[\xymatrix{
D_1\otimes E_1\ar[rr]^{\tau_D\otimes 1}\ar[d]_{c'_{D_1, E_1}}&&D_2\otimes E_1\ar[rr]^{1\otimes\tau_E}\ar@<-1ex>[d]_{c'_{D_2, E_1}}&&D_2\otimes E_2\ar@<1ex>[d]^{c'_{D_2, E_2}}\\
E_1\otimes D_1\ar[rr]_{1\otimes\tau_D}&&E_1\otimes D_2\ar[rr]_{\tau_E\otimes 1}\ar@<-1ex>[u]_{c'_{E_1, D_2}}&&E_2\otimes D_2\ar@<1ex>[u]^{c'_{E_2, D_2}}
}\]
The pair of arrows on the central vertical trivially compose to identity, and the left, right (inner) rectangular subdiagrams are the squares of naturality for $c'$. Composition of the pair of vertical arrows on the far right is $1\tens{}c^2_{E, D}c^2_{D, E}$, and this is identity because we identify the symmetric tensor category $\mc{D}$ with its image under the given braided inclusion to $\mc{Z}(\mc{C}_2)$. We could have written the horizontal compositions in the opposite sense, in which case we would have needed $\mc{D}$ symmetric tensor subcategory of $\mc{Z}(\mc{C}_1)$. Thus $\mc{D}$ sits inside $\mc{Z}(\mc{C}_1\tens{}\mc{C}_2)$ canonically (up to braided equivalence) whenever both $\mc{C}_i$ are tensor over $\mc{D}$. 
%
\end{proof}
\begin{remark}\label{EquivTensProductCenter}
The canonical tensor category structure on $\mc{C}_1\tens{D}\mc{C}_2$ which is the content of Corollary \ref{monoidalityOfTensor} may be described explicitly as follows. Write $\otimes_i$ for the tensor product and $c^i$ for the braiding in $\mc{C}_i$, respectively. Denote by $B_{1, 2}$ the universal $\mc{D}$-balanced functor from $\mc{C}_1\tens{}\mc{C}_2$. Let $\hat{\tau}:(\mc{C}_1\tens{}\mc{C}_2)^{\tens{}2}\rightarrow(\mc{C}_1\tens{}\mc{C}_2)^{\tens{}2}$ be the equivalence of the form $(13)(24)$, i.e. $\hat{X}\tens{}\hat{Y}\mapsto\hat{Y}\tens{}\hat{X}$ for $\hat{X} = X_1\tens{}X_2\in\mc{C}_1\tens{}\mc{C}_2$, and let $\tau:\mc{C}_1\tens{}\mc{C}_2\rightarrow\mc{C}_2\tens{}\mc{C}_1$ be the switch $(12)$. Let $\otimes = \otimes_1\tens{}\otimes_2$ be the tensor product in $\mc{C}_1\tens{}\mc{C}_2$. Being completely explicit, the $\mc{D}$-balancing of $B_{1,2}\otimes:(\mc{C}_1\tens{}\mc{C}_2)^{\tens{}2}\rightarrow\mc{C}_1\tens{D}\mc{C}_2$ in positions 1 and 3 looks like
\begin{eqnarray*}
b_{X_1Y_1, D, X_2Y_2}\circ B_{1,2}(1\otimes c^1_{D, Y_1}\tens{}1):&&\!\!\!\!\!\!\!\!\!\!B_{1,2}(X_1DY_1\tens{}X_2Y_2)\rightarrow B_{1,2}(X_1Y_1\tens{}DX_2Y_2)\\
B_{1,2}(1\tens{} c^2_{D, X_2}\otimes 1)\circ b_{X_1Y_1, D, X_2Y_2}:&&\!\!\!\!\!\!\!\!\!\!B_{1,2}(X_1Y_1D\tens{}X_2Y_2)\rightarrow B_{1,2}(X_1Y_1\tens{}X_2DY_2)
\end{eqnarray*}
respectively. Now abbreviate $\hat{X}D:=(X_1\otimes D)\tens{}X_2$ and $D\hat{X}: = X_1\tens{}(D\otimes X_2)$ for any object $D\in\mc{D}$. Balancing of $B_{1,2}\otimes$ in positions 1 and 3 can then be written as the following compositions on the level of objects in $\mc{C}_1\tens{}\mc{C}_2$:
\begin{eqnarray}
&&B_{1,2}\otimes(\hat{X}D\tens{}\hat{Y})\rightarrow B_{1,2}\otimes(\hat{X}\tens{}\hat{Y}D)\rightarrow B_{1,2}\otimes(D\hat{X}\tens{}\hat{Y})\label{1}\\
&&B_{1,2}\otimes(\hat{X}\tens{}\hat{Y}D)\rightarrow B_{1,2}\otimes(D\hat{X}\tens{}\hat{Y})\rightarrow B_{1,2}\otimes(\hat{X}\tens{}D\hat{Y})\label{2}.
\end{eqnarray}
Written this way it's easy to see how the components of the various balancing structures act. In order to see that these balancings are consistent in the sense of Equation (\ref{MultiBalancedCondition2}) consider the diagram below. To save space we have suppressed $B_{1,2}\otimes$ at each vertex and have adopted the convention of suppressing components of morphisms which are irrelevant to the commutativity of the diagrams of which they form an edge.
\[\!\!\!\!\!\!\!\!\!\!\!\!\!\!\!\!\!
\xymatrix{
E\hat{X}D\tens{}\hat{Y}\ar[drrr]_{1_{X_1}\otimes c^1_{D, Y_1}}\ar[rrr]^{c^2_{E, X_2}\otimes 1_{Y_2}}&&&\hat{X}D\tens{}E\hat{Y}\ar[rrr]^{1_{X_1}\otimes c^1_{D, Y_1}}&&&\hat{X}\tens{}E\hat{Y}D\ar[d]^{b_D}\\
\hat{X}D\tens{}\hat{Y}E\ar[drr]|{1_{X_1}\otimes c^1_{D, Y_1}\otimes 1_E}\ar[dd]|<<<<<<<<<{1_{X_1}\otimes c^1_{D, Y_1E}}\ar[u]_{b_E}&&&E\hat{X}\tens{}\hat{Y}D\ar[urrr]_{c^2_{E, X_2}\otimes 1_{Y_2}}\ar[dr]^{b_D}&&&D\hat{X}\tens{}E\hat{Y}\\
&&\hat{X}\tens{}\hat{Y}DE\ar[ur]^{b_{E}}\ar@<.5ex>[rr]^{b_{DE}}\ar@<-.5ex>[rr]_{b_{ED}}\ar@{.>}@<-.5ex>[dll]_>>>>>>>>>>>{1_{X_1}\otimes 1_{Y_1}\otimes c^1_{D, E}\quad}&&DE\hat{X}\tens{}\hat{Y}\ar[urr]|{1_D\otimes c^2_{E, X_2}\otimes 1_{Y_2}}\ar@{.>}@<-.5ex>[drr]_>>>>>>>>>>>>{c^2_{D, E}\otimes 1_{X_2}\otimes 1_{Y_2}}
&&\\
\hat{X}\tens{}\hat{Y}ED\ar[rrr]_{b_D}\ar@{.>}@<-.5ex>[urr]_{\quad1_{X_1}\otimes 1_{Y_1}\otimes c^1_{E, D}}&&&D\hat{X}\tens{}\hat{Y}E\ar[rrr]_{b_E}&&&ED\hat{X}\tens{}\hat{Y}\ar[uu]|>>>>>>>>>{c^2_{E, DX_2}\otimes 1_{Y_2}}\ar@{.>}@<-.5ex>[ull]_<<<<<<<<<<{\,\,\,\,\,\,c^2_{E, D}\otimes 1_{X_2}\otimes 1_{Y_2}}
}
\]
Each edge labeled by a component of $c^i$ corresponds to the central structure of the image of $\mc{D}$ in $\mc{Z}(\mc{C})$. Every subdiagram is either naturality of the balancing $b$, the balancing diagram in Definition \ref{DEF;bal}, or commutes trivially. Each double perforated-arrow edge composes to identity (see the proof of Corollary \ref{monoidalityOfTensor} above). The perimeter is the analogue of Equation (\ref{MultiBalancedCondition2}) for $n=4$, $i=1$ and $j=3$. Thus $B_{1,2}\otimes$ is balanced in positions 1 and 3 simultaneously and thus descends to the functor $\overline{\otimes}:(\mc{C}_1\tens{D}\mc{C}_2)^{\tens{}2}\rightarrow\mc{C}_1\tens{D}\mc{C}_2$ satisfying the properties described in Definition \ref{DEF;tensor}. This discussion is essentially contained in Section 6 of \cite{JG}.
\end{remark}
\begin{remark}
Suppose $\mc{D}$ is braided with braiding $d$, and let $A, B$ be two commutative algebras in $\mc{D}$. Let $\mc{D}_A$ be the category of left $A$-modules in $\mc{D}$, and similarly define $\mc{D}_B$. Then $\mc{D}_A$, $\mc{D}_B$ are $\mc{D}$-bimodule categories, and we have 
\begin{equation*}
\mc{D}_A\tens{D}\mc{D}_B\simeq\mc{D}_{A\otimes B}.
\end{equation*} 
In order for $A\otimes B$ to be itself a commutative algebra it is necessary that $A$ and $B$ centralize each other in $\mc{D}$, i.e. $c_{A, B}=c^{-1}_{B, A}$ (see Section \ref{relCentBraidSect} for further discussion). In order for the category of modules $\mc{D}_{A\otimes B}$ to have the structure of a monoidal category it is sufficient that $A\in\mc{D}'$. Since this discussion is symmetric in $A, B$ (due to the braiding of $\mc{D})$ we must also have $B\in\mc{D}'$. We restate this observation in the next proposition. The author thanks Alexei Davydov for bringing his attention to this result in \cite{WittBraided}.
\begin{prop}
Let $\mc{D}$ be a braided monoidal category, and suppose that $A, B$ are commutative algebras in $\mc{D}$. Then the category of modules $\mc{D}_{A\otimes B}$ has the structure of a monoidal category over $\mc{D}$ iff both $A$ and $B$ are in $\mc{D}'$. 
\end{prop}
Of  course it would be fruitful to discuss general categorical properties reflected by the mutual centralization of the algebras $A$ and $B$. Unfortunately we must defer this discussion to a later article. We refer the reader to \textit{loc. cit.} where related issues are addressed. 
\end{remark}
%
\section{Braiding the relative tensor product}\label{relCentBraidSect}
The following definition first appeared in papers of M\"{u}ger (\cite{Mu1}, \cite{Mu2}) and was used to great effect in \cite{braidedI}. Let $\mc{C}$ be a braided tensor category with braiding $c$. Objects $X, Y$ in $\mc{C}$ are said to \textit{centralize} each other if 
\begin{equation}\label{centralizer equation}
c_{Y, X}\circ c_{X, Y}=id_{X\otimes Y}.
\end{equation}
\begin{defn}\label{centralizerDefn}
Let $\mc{C}$ be a braided tensor category. For tensor subcategory $\mc{D}\subset\mc{C}$ the \textit{centralizer} of $\mc{D}$ in $\mc{C}$, denoted $\mc{D}'$,  is defined to be the full subcategory of $\mc{C}$ consisting of objects centralizing every object in $\mc{D}$. 
\end{defn}
Evidently $\mc{D}'$ is fusion whenever $\mc{C}$ is fusion and $\mc{D}$ is a fusion subcategory. In the case $\mc{D}=\mc{C}$ the centralizer $\mc{C}'$ is sometimes referred to as the M\"uger center of $\mc{C}$. We will call it simply \textit{the} centralizer of $\mc{C}$. Recall Definition \ref{braidingOVERdef}: for $\mc{C}, \mc{D}$ braided tensor categories we say $\mc{C}$ is braided \textit{over} $\mc{D}$ if there is a braided inclusion $\mc{D}\hookrightarrow\mc{C}'$. 
\begin{remark}[\cite{braidedI}]\label{BraidedOverIsTensorOver} Any braided tensor category $\mc{C}$ over $\mc{D}$ as in Definition \ref{braidingOVERdef} can also be viewed as a tensor category over $\mc{D}$ as follows. Define the braided inclusion $\mc{D}\hookrightarrow\mc{Z}(\mc{C})$ by the composition $\mc{D}\hookrightarrow\mc{C}'\hookrightarrow\mc{C}\rightarrow\mc{Z}(\mc{C})$ where the last step (which is fully faithful) comes from the braiding in $\mc{C}$. 
\end{remark}
Continuing with the notational conventions introduced in Remark \ref{EquivTensProductCenter} we are ready to address our attention to the proof of Theorem \ref{braiding of relative center thm}.
\begin{proof}[Proof of Theorem \ref{braiding of relative center thm}] 
By hypothesis we have braided inclusions $\mc{D}\hookrightarrow\mc{C}_i$. By Remark \ref{BraidedOverIsTensorOver} we can view the braided tensor categories $\mc{C}_1$, $\mc{C}_2$ as being each tensor over $\mc{D}$. As discussed in the proof of Corollary \ref{monoidalityOfTensor} this allows us to view $\mc{C}_1\tens{}\mc{C}_2$ as being tensor over $\mc{D}$, and thus the discussion about the tensor category structure on the relative tensor product $\mc{C}_1\tens{D}\mc{C}_2$, either there or in Remark \ref{EquivTensProductCenter}, applies here. 

Denote by $\hat{\tau}:(\mc{C}_1\tens{}\mc{C}_2)^{\tens{}2}\rightarrow(\mc{C}_1\tens{}\mc{C}_2)^{\tens{}2}$ the functor which acts by permuting the factors of the Deligne product by the action of the permutation $(13)(24)$. One shows that $B_{1,2}\otimes\hat{\tau}$ is balanced in positions 1 and 3 simultaneously in precisely  the same way we showed this for $B_{1,2}\otimes$ in Remark \ref{EquivTensProductCenter}. 
It therefore descends to the functor $\overline{\otimes}^{rev}$. Now write $c^i$ for the braiding in $\mc{C}_i$ and denote by $c'$ the braiding in $\mc{C}_1\tens{}\mc{C}_2$ given by
\begin{equation}\label{TwistBraid}
c'_{\hat{X}, \hat{Y}} = c^1_{X_1, X_2}\tens{}c^2_{Y_1, Y_2}:\hat{X}\otimes\hat{Y}\stackrel{\sim}{\rightarrow}\hat{Y}\otimes\hat{X}.
\end{equation}
This allows us to define the natural isomorphism $B_{1,2}*c':B_{1,2}\otimes\stackrel{\sim}{\rightarrow}B_{1,2}\otimes^{rev}$ having components $(B_{1,2}*c')_{\hat{X}, \hat{Y}} = B_{1,2}(c^1_{X_1, X_2}\tens{}c^2_{Y_1, Y_2})$. We show that it is $\mc{D}$-balanced in positions 1 and 3 (recall Remark 
\ref{multi-relative tensor}). 
The diagram below shows that $B_{1,2}*c'$ is $\mc{D}$-balanced in position 1. 
\[
\xymatrix{
B_{1,2}\otimes(\hat{X}D\tens{}\hat{Y})\ar[d]_{B_{1,2}(1\otimes c^1_{D, Y_1}\tens{}1)}\ar[rrrrr]^{(B_{1,2}*c')_{\hat{X}D\tens{}\hat{Y}}}&&&&&B_{1,2}\otimes(\hat{Y}\tens{}\hat{X}D)\ar[d]^{b_{Y_1X_1, D, Y_2X_2}}\\
B_{1, 2}\otimes(\hat{X}\tens{}\hat{Y}D)\ar[d]_{b_{X_1Y_1, D, X_2Y_2}}\ar[urrrrr]|{B_{1,2}(c^1_{X_1, Y_1}\otimes 1\tens{}1)}&&&&&B_{1, 2}\otimes(D\hat{Y}\tens{}\hat{X})\ar[d]^{B_{1,2}(1\tens{}c^2_{D, Y_2}\otimes 1)}\\
B_{1,2}\otimes(D\hat{X}\tens{}\hat{Y})\ar[urrrrr]|{B_{1,2}(1\tens{}1\otimes c^2_{X_2, Y_2})}\ar[rrrrr]_{(B_{1,2}*c')_{D\hat{X}\tens{}\hat{Y}}}&&&&&B_{1,2}(\hat{Y}\tens{}D\hat{X})
}
\]
An analogous diagram shows that $B_{1,2}*c'$ is $\mc{D}$-balanced in the third position. Thus we have a unique natural isomorphism $\overline{c}:=\overline{B_{1,2}*c'}$ satisfying $\overline{c}*(B_{1, 2}\tens{}B_{1, 2})=B_{1, 2}*c$ (recall Definition \ref{DescendedBalancedNatTrans} and the following Remark \ref{multi-relative tensor}) as in the diagram of 2-cells below:
\[
\xymatrix{
(\mc{C}_1\tens{}\mc{C}_2)^{\tens{}2}\ar[dd]_{B_{1,2}\tens{}B_{1,2}}\ar[rrrr]^{\hat{\tau}}\ar@{~>}[ddrrrr]^>>>>>>>>>>>>>>>{}="3"_\otimes&&&&(\mc{C}_1\tens{}\mc{C}_2)^{\tens{}2}\ar[dd]_>>>>>>>>>>>>{}="4"^\otimes\\
&&&&\\
(\mc{C}_1\tens{D}\mc{C}_2)^{\tens{}2}\ar@/_1pc/[rr]_{\overline{\otimes}^{rev}}^{}="2"\ar@{~>}@/^1pc/[rr]^{\overline{\otimes}}_{}="1"&&\mc{C}_1\tens{D}\mc{C}_2&&\mc{C}_1\tens{}\mc{C}_2\ar@<.5ex>[ll]^{B_{1,2}}\ar@<-.5ex>@{~>}[ll]
\ar@{=>}"1";"2"^{\,\,\overline{c}}
\ar@/^.7pc/@{=>}"3";"4"^{c'}
}
\]
The squiggly arrows form a diagram together with the left-most vertical, and the remaining bland arrows form a separate diagram. These two diagrams are the definitions of $\overline{\otimes}$ and $\overline{\otimes}^{rev}$, respectively, and commute on their own without recourse to the indicated 2-cells. The natural isomorphism $\overline{c}$ constitutes a braiding in $\mc{C}_1\tens{D}\mc{C}_2$; braiding diagrams are satisfied by $\overline{c}$ because they are satisfied by $c$. 
\end{proof}
\begin{remark}
Theorem \ref{braiding of relative center thm} also follows from Corollary \ref{monoidalityOfTensor} and a result of Joyal and Street in \cite{JS}. We briefly recall the relation between braiding and multiplication structures on monoidal categories discussed there. Let $\mc{C}$ be a monoidal category. Then the cartesian product $\mc{C}\times\mc{C}$ has the structure of a monoidal category with monoidal structure
\[
\xymatrix{
\mc{C}\times\mc{C}\times\mc{C}\times\mc{C}\ar[rr]^{id\times\textrm{\,switch\,}\times id}&&\mc{C}\times\mc{C}\times\mc{C}\times\mc{C}\ar[r]^<<<<<{\otimes\times\otimes}&\mc{C}\times\mc{C}
}\]
which we will here denote $\otimes'$. A \textit{multiplication} on $\mc{C}$ is a monoidal functor $\mc{C}\times\mc{C}\rightarrow\mc{C}$ together with natural unit cohesion isomorphisms required to satisfy the expected diagrams. Details may be found in the original paper. 

The observation important for us here is the correspondence between multiplications on $\mc{C}$ and braidings in $\mc{C}$ (Propositions 5.2 and 5.3, \textit{loc. cit.}). One may show that the tensor structure $\overline{\otimes}:\mc{C}_1\tens{D}\mc{C}_2\times\mc{C}_1\tens{D}\mc{C}_2\rightarrow\mc{C}_1\tens{D}\mc{C}_2$ is in fact a multiplication in $\mc{C}_1\tens{D}\mc{C}_2$ under the stated hypotheses with monoidal structure coming from the braidings in $\mc{C}_i$. By \cite{JS} it then corresponds to a braiding, and this braiding is precisely the one described in the proof of Theorem  \ref{braiding of relative center thm}. 
\end{remark}
\section{Examples}\label{ExampleSection} The most basic examples of the phenomenon described in the previous sections occur for braided categories structurally determined by various group theoretical data. In this section we consider braided fusion categories related to pre-metric groups (recall \S \ref{PreMetGrpSec} and associated citations) and categories of representations of finite groups. In order to make the necessary calculations we will need a lemma on tensor functors. We state the result without proof as it is readily located in the literature or else easily proved.
\begin{lem}\label{SurjTensFunctLem}
Let $\mc{C}, \mc{D}$ be finite tensor categories, and let $F:\mc{C}\rightarrow\mc{D}$ be a surjective tensor functor having right adjoint $I$. Then $\FPdim(\mc{C})=\FPdim(\mc{D})\FPdim(I(1))$.
\end{lem}
\begin{ex}[Pre-metric groups]\label{PreMetEx}
Recall \S \ref{PreMetGrpSec}. Let $G_i$, $i=1, 2$, be finite abelian groups, and let $q_i:G_i\rightarrow k^\times$ be quadratic forms satisfying $q_i(g)=\beta_i(g, g)$ for some bicharacters $\beta_i$ on $G_i$. 
One easily checks that $(G_1\times G_2, p)$ is a pre-metric group for $p(g, h)=q_1(g)q_2(h)$. As a quadratic form $p$ comes from the bicharacter on $G_1\times G_2\times G_1\times G_2$ given by $(g_1, h_1, g_2, h_2)\mapsto\beta_1(g_1, g_2)\beta_2(h_1, h_2)$. Likewise, if $b_i$ are the symmetric bimultiplicative forms determined by $q_i$ then $b((g_1, h_1), (g_2, h_2)):=b_1(g_1, g_2)b_2(h_1, h_2)$ is that determined by $p$. 

Now suppose we have embeddings $G\hookrightarrow G_i$ for a finite group $G$. Supposer further that these embeddings satisfy $q_1(g)=q_2(g)$ for all $g\in G$. Then the pair $(G, q)$ is a metric group for $q:=q_i|_{G}$. Denote by $\tilde{G}$ the subgroup of $G_1\times G_2$ given by the set $\{(x, x^{-1})|x\in G\}$ and suppose that $p$ descends to a quadratic form on $(G_1\times G_2)/\tilde{G}$. Since $p$ is constant on $\tilde{G}$-cosets we have $q_i(x)q(g)
=q_i(gx)q(1)$ for $x\in G_i, g\in G$ as the reader can easily verify. As a result $b_i(g, x)=q(1)^{-1}$, and since $b_i(1, 1)=1$ we conclude that
\begin{equation}\label{centerPreMetGrp}
b_i(g, x)=1
\end{equation}
for each $i=1, 2$ whenever $x\in G_i$ and for $g$ any element of $G$.

Let's translate this into the language of pointed braided fusion categories via \S \ref{PreMetGrpSec}. Pre-metric inclusions $G\hookrightarrow G_i$ correspond to braided inclusions $\mc{C}(G, q)\rightarrow\mc{C}(G_i, q_i)$. Equation (\ref{centerPreMetGrp}) becomes $c^i_{X, Y}=1$ ($c^i$ the braiding in $\mc{C}(G_i, q_i)$) whenever $X, Y$ are homogeneous objects of $Vec_{G_i}$ of degrees $g\in G$, $x$ respectively. Thus the images of the braided inclusions are contained in M\"uger centers $\mc{C}(G_i, q_i)'$, i.e. $\mc{C}(G_i, q_i)$ are braided over $\mc{C}(G, q)$. 

We will now show that $\mc{C}(G_1\times G_2/\tilde{G}, p)\simeq \mc{C}(G_1, q_1)\boxtimes_{\mc{C}(G, q)}\mc{C}(G_2, q_2)$ as braided fusion categories. Note that $\mc{C}(G_1, q_1)$ is a right $\mc{C}(G, q)$-module category; for $g\in G_1$ and $h\in G$ define $(V_g, W_h)\mapsto (V\otimes W)_{gh}$ where $X_y$ denotes a homogeneous space $X$ of degree $y$. Similarly $\mc{C}(G_2, q_2)$ is a left $\mc{C}(G, q)$-module category. Denote by $\overline{(x, y)}$ the coset represented by $(x, y)$ in the factor group $G_1\times G_2/\tilde{G}$. The braiding on $\mc{C}(G_1\times G_2/\tilde{G}, p)$ is given by $c_{X_{\overline{(g, h)}}, V_{\overline{(e, f)}}}=\beta^1(g, e)\beta^2(h, f)\tau_{X, V}$ so the $\mc{C}(G, q)$-balanced functor 
\begin{equation*}
F:\mc{C}(G_1, q_1)\tens{}\mc{C}(G_2, q_2)\rightarrow\mc{C}(G_1\times G_2/\tilde{G}, p),\quad\quad V_g\tens{}W_h\mapsto(V\otimes W)_{\overline{(g, h)}}
\end{equation*}
is braided. The functor $F$ descends to a unique surjective braided functor 
\begin{equation*}
\overline{F}: \mc{C}(G_1, q_1)\boxtimes_{\mc{C}(G, q)}\mc{C}(G_2, q_2)\rightarrow\mc{C}(G_1\times G_2/\tilde{G}, p)
\end{equation*}
via the surjective universal $\mc{C}(G, q)$-balanced tensor functor $B:\mc{C}(G_1, q_1)\boxtimes_{}\mc{C}(G_2, q_2)\rightarrow\mc{C}(G_1, q_1)\boxtimes_{\mc{C}(G, q)}\mc{C}(G_2, q_2)$ (Definition \ref{DEF;tensor}). We show that $\overline{F}$ is actually an equivalence. Let $I$ be right adjoint to $B$. It is straightforward to check that
\begin{equation*}
IB(1\tens{}1)=\bigoplus_{g\in G}k_g\tens{}k_{g^{-1}}
\end{equation*}
and consequently $\FPdim(IB(1\tens{}1))=|G|$. By Lemma \ref{SurjTensFunctLem} we therefore have the equation
\begin{equation*}
\FPdim(\mc{C}(G_1, q_1)\boxtimes_{\mc{C}(G, q)}\mc{C}(G_2, q_2))=\frac{\FPdim(\mc{C}(G_1, q_1))\FPdim(\mc{C}(G_2, q_2))}{\FPdim(\mc{C}(G, q))}.
\end{equation*}  
The expression on the right is $|G_1||G_2|/|G|=\FPdim(\mc{C}(G_1\times G_2/\tilde{G}, p))$ and thus $\overline{F}$ is an equivalence. 
\end{ex}
\begin{ex}[Fibre product of groups]Let $G_1, G_2$ be finite groups having normal subgroups $N_1, N_2$, respectively, with isomorphic quotient groups ($G_1/N_1\simeq G_2/N_2$). Let $\phi:G:=G_2/N_2\rightarrow G_1/N_1$ be an isomorphism and let $\pi_i:G_i\rightarrow G_i/N_i$ be the canonical projections. Then we may form the fibre product $G_1\times_G G_2$ as the pullback of $\pi_1$ and $\phi\circ\pi_2$ in the usual way (rectangular subdiagram below).
\[\xymatrix{
G_1\times_G G_2\ar[rr]^{\theta_1}\ar[d]_{\theta_2}&&G_1\ar@{->>}[d]^{\pi_1}&\\
G_2\ar@{->>}[r]_{\pi_2}&G_2/N_2\ar[r]^\sim_{\varphi}&G_1/N_1\ar[r]^{\kappa}&\End(V)
}\]
Let $\iota:G_1\times_GG_2\hookrightarrow G_1\times G_2$ be the inclusion of the fibre product as a subgroup of the cartesian product of groups. Then we have a (canonical) surjective monoidal functor $\overline{\iota}:\Rep(G_1\times G_2)\rightarrow\Rep(G_1\times_G G_2)$ which acts by $\iota$-precomposition. The category $\Rep(G)$ acts on $\Rep(G_1)$ from the right and on $\Rep(G_2)$ from the left by the formulas
\begin{equation*}
\rho_1\triangleleft\kappa:=\rho_1\otimes\kappa\pi_1,\qquad \kappa\triangleright\rho_2:=\kappa\varphi\pi_2\otimes\rho_2
\end{equation*}
where $\rho_i$ is a representation of $G_i$ and $\kappa$ any representation of $G_1/N_1$. We show that $\overline{\iota}$ is $\Rep(G)$-balanced as a functor from $\Rep(G_1)\tens{}\Rep(G_2)$. For any pair $(g, h)\in G_1\times_G G_2$ we have $\kappa\pi_1(g)=\kappa\varphi\pi_2(h)$ by definition, hence
\begin{eqnarray*}
\overline{\iota}(\rho_1\triangleleft\kappa, \rho_2)(g, h)&=&\rho_1(g)\otimes\kappa\pi_1(g)\otimes\rho_2(h)\\
&=&\rho_1(g)\otimes\kappa\varphi\pi_2(h)\otimes\rho_2(h)\\
&=&\overline{\iota}(\rho_1, \kappa\triangleright\rho_2)(g, h).
\end{eqnarray*}
That is, $\overline{\iota}:\Rep(G_1)\tens{}\Rep(G_2)\rightarrow\Rep(G_1\times_GG_2)$ is $\Rep(G)$-balanced (the balancing constraint is trivial). The unique extension $\overline{\overline{\iota}}$ in the diagram below is also surjective.
\[
\xymatrix{
\Rep(G_1)\tens{}\Rep(G_2)\ar[dr]^{\overline{\iota}}\ar[d]_B&\\
\Rep(G_1)\boxtimes_{\Rep(G)}\Rep(G_2)\ar[r]_>>>>>{\overline{\overline{\iota}}}&\Rep(G_1\times_GG_2)
}\]
To show that $\overline{\overline{\iota}}$ is an equivalence it suffices to show that the Frobenius-Perron dimensions of the domain and codomain categories are equal (this is the same as the treatment of $\overline{F}$ in Example \ref{PreMetEx} above and works here for the same reason, i.e. that $\overline{\overline{\iota}}$ is a tensor functor). We have
\begin{equation*}
\FPdim(\Rep(G_1\times_GG_2))=\displaystyle{\frac{|G_1||G_2|}{|G|}}=\displaystyle{\frac{\FPdim(\Rep(G_1)\boxtimes\Rep(G_2))}{\FPdim(\Rep(G))}}.
\end{equation*}
Since $\FPdim(I(1))=\FPdim(\Rep(G))=|G|$ for $I$ the right adjoint of the universal surjective $\Rep(G)$-balanced tensor functor $B$ the expression at far right is precisely $\FPdim(\Rep(G_1)\boxtimes_{\Rep(G)}\Rep(G_2))$. 
\end{ex}
%
\section{Module categories over the relative center}\label{LemmasSection}
In this section we study module categories over tensor categories of the form $\cent{C}{D}$. We begin with a general lemma relating balancing and module category structure. Let $\mc{D}$ be a symmetric tensor category. 
In the case that $\mc{C}$ is tensor over $\mc{D}$ it is natural to ask about the relation between $\mc{C}$-module and $\cent{C}{D}$-module categories. This is the subject of the first proposition in this section.
\begin{prop}\label{FunctorialityOfZ_DLemma}
Let $\mc{C}$ be tensor over $\mc{D}$, and suppose $\mc{M}$ is a $\mc{C}$-bimodule category. Then $\cent{M}{D}$ has canonical structure of a $\cent{C}{D}$-bimodule category. If $G:\mc{M}\rightarrow\mc{N}$ is a $\mc{C}$-module functor then there is a canonical $\cent{C}{D}$-bimodule functor $\mc{Z}_{\mc{D}}(G):\cent{M}{D}\rightarrow\cent{N}{D}$.
\end{prop}
\begin{proof} The module category $\mc{M}$ is a $\mc{D}$-bimodule category by restriction of the $\mc{C}$-module category structure to the image of $\mc{D}$ in $\mc{C}$. Define an action of $\cent{C}{D}$ on $\cent{M}{D}$ by 
\begin{equation}\label{centralactiondef}
(X, c_X)\otimes(M, \varphi_M):=(X\otimes M, c_X*\varphi_M)
\end{equation}
where for any $Y\in\mc{D}$ the isomorphisms $c_X*\varphi_M$ are defined by the composition 
\begin{equation}\label{centralactiondefII}
(c_X*\varphi_M)_Y:=(id_X\otimes\varphi_{Y, M})(c_{Y, X}\otimes id_M): Y\otimes(X\otimes M)\rightarrow(X\otimes M)\otimes Y
\end{equation}
whenever $Y\in\mc{D}$ and $X\in\mc{C}$. The diagram below right shows that $c_X*\varphi_M$ satisfies the relative braiding condition (the diagram in Definition \ref{defCent}).
\[\!\!\!\!\!\!\!\!\!\!\!\!\!\!\!\!\!\!\!\!\xymatrix{
YWXM\ar[rr]^{(c_{WX}*\varphi_M)_Y}\ar[dd]_{c_{Y, W}}&&WXMY\\
&&\\
WYXM\ar[rr]_{c_{Y, X}}&&WXYM\ar[uu]_{\varphi_{Y, M}}
}\;\quad
\xymatrix{
YZXM\ar[rr]^{c_{YZ, X}}\ar[dr]^{c_{Z, X}}\ar[dd]_{(c_X*\varphi_M)_Z}&&XYZM\ar[dr]^{\varphi_{YZ, M}}\ar[dd]_{\varphi_{Z, M}}&\\
&YXZM\ar[ur]^{c_{Y, X}}\ar[dl]^{\varphi_{Z, M}}&&XMYZ\\
YXMZ\ar[rr]_{c_{Y, X}}&&XYMZ\ar[ur]_{\varphi_{Y, M}}&
}\]
Every subdiagram on right is either one of the braiding diagram in $\cent{C}{D}$ or $\cent{M}{D}$, or else commutes trivially. The path along the top is the definition of $(c_X*\varphi_M)_{Y\otimes Z}$ and that along the bottom is $(c_X*\varphi_M)_Y\otimes id_Z$. The diagram on left shows that the action is strictly associative: the three-arrow composition from upper left entry to upper right entry is precisely $(c_W*(c_X*\varphi_M))_Y$ for $W, X\in\mc{C}$ and $Y\in\mc{D}$. Hence $\cent{M}{D}$ has the structure of a (strict) left $\cent{C}{D}$-module category (strictness of associativity also follows from Remark \ref{2-remark} below). Verification that the right action 
\begin{equation*}\label{centralCenterAction}
(M, \varphi_M)\otimes (X, c_X):=(M\otimes D, \varphi_M*c_X),\quad (\varphi_M*c_X)_Y:=(id_M\otimes c_{Y, X})(\varphi_{Y, M}\otimes id_X)
\end{equation*}
gives $\cent{M}{D}$ the structure of a right $\cent{C}{D}$-module category is identical. The left and right module structures are trivially consistent, hence combine to give $\cent{M}{D}$ the structure of a $\cent{C}{D}$-bimodule category.

Finally let $G:\mc{M}\rightarrow\mc{N}$ be a $\mc{C}$-bimodule functor with bimodule structure $g$. Then define functor $\mc{Z}_\mc{D}(G):\cent{M}{D}\rightarrow\cent{N}{D}$ by $(M, \varphi_M)\mapsto(G(M), \varphi^G_{M})$ where the components of $\varphi^G_{M}$ are defined by the composition
\[\xymatrix{
\varphi^G_{Y, M}:=Y\otimes G(M)\ar[r]^>>>>>{g^{-1}_{Y\tens{}1, M}}&G(Y\otimes M)\ar[r]^{G(\varphi_{Y, M})}&G(M\otimes Y)\ar[r]^{g_{1\tens{}Y, M}}&G(M)\otimes Y.
}\]
The following diagram shows that $\varphi^G_{M}$ satisfies the braiding axiom (Definition \ref{defCent}). Let $Y, Z\in\mc{D}$.
\[\!\!\!\xymatrix{
YZG(M)\ar[r]^{g^{-1}_{YZ\tens{}1, M}}\ar[d]_{g^{-1}_{Y\tens{}1, M}}&G(YZM)\ar[rr]^{G(\varphi_{YZ, M})}\ar[dr]_{G(1\otimes\varphi_{Z, M})}\ar[dl]^{g_{Y\tens{}1, ZM}}&&G(MYZ)\ar[r]^{g_{1\tens{}YZ, M}}\ar[dr]_{g_{1\tens{}Z, MY}}&G(M)YZ\\
YG(ZM)\ar[d]_{G(\varphi_{Z, M})}&&G(YMZ)\ar[dll]_{g_{Y\tens{}1, MZ}}\ar[d]_{g_{Y\tens{}Z, M}}\ar[ur]_{G(\varphi_{Y, M}\otimes 1)}\ar[drr]^{g_{1\tens{}Z, YM}}&&G(YM)Z\ar[u]_{g_{1\tens{}Y, M}}\\
YG(MZ)\ar[rr]_{g_{1\tens{}Z, M}}&&YG(M)Z\ar[rr]_{g^{-1}_{Y\tens{}1, M}}&&G(YM)Z\ar[u]_{G(\varphi_{Y, M})}
}\]
Every subdiagram is either $\mc{C}$-bimodule linearity for $G$ (restricted to $\mc{D}$), naturality of $g$ or the braiding diagram for $\varphi$ in $\cent{M}{D}$. The perimeter is the equation $\varphi_{YZ, G(M)}=(\varphi_{Y, G(M)}\otimes 1_Z)(1_Y\otimes\varphi_{Z, G(M)})$.

To see that $\mc{Z}_{\mc{D}}(G)$ is a $\cent{C}{D}$-module functor consider the following diagram where $(X, c_X)\in\cent{C}{D}$ and $Y\in\mc{D}$.
\[\xymatrix{
YG(XM)\ar[r]^{g^{-1}_{Y\tens{}1, XM}}\ar[d]_{g_{X\tens{}1, M}}&G(YXM)\ar[dl]|{g_{YX, M}}\ar[r]^{G(c_{Y, X})}&G(XYM)\ar[r]^{G(\varphi_{Y, M})}\ar[dl]|{g_{XY, M}}\ar[d]|{g_{X\tens{}1, YM}}&G(XMY)\ar[d]|{g_{X\tens{}1, MY}}\ar[r]^{g_{1\tens{}Y, XM}}\ar[dr]|{g_{X\tens{}Y, M}}&G(XM)Y\ar[d]^{g_{X\tens{}1, M}}\\
YXG(M)\ar[r]_{c_{Y, X}}&XYG(M)\ar[r]_{g^{-1}_{Y\tens{}1, M}}&XG(YM)\ar[r]_{G(\varphi_{Y, M})}&XG(MY)\ar[r]_{g_{1\tens{}Y, M}}&XG(M)Y
}\]

Every subdiagram is either naturality or the bimodule condition on $g$. The perimeter is the equation
\begin{equation*}
(g_{X\tens{}1, M}\otimes 1_Y)(c_X*\varphi_M)^G_Y=(c_X*\varphi^G_M)_Y(1_Y\otimes g_{X\tens{}1, M}).
\end{equation*}
Thus the components $g_{X\tens{}1, M}:(G(X\otimes M), (c_X*\varphi_M)^G)\rightarrow (X\otimes G(M), c_X*\varphi_M^G)$ form a family of natural isomorphisms 
\begin{equation*}
\mc{Z}_\mc{D}(G)((X, c_X)\otimes (M, \varphi_M))\stackrel{\sim}{\rightarrow}(X, c_X)\otimes\mc{Z}_\mc{D}(G)(M, \varphi)
\end{equation*}
in $\cent{N}{D}$ and constitute a left $\cent{C}{D}$-module structure for the functor $\mc{Z}_\mc{D}(G)$. The necessary diagrams come from those satisfied by $g$ as $\mc{C}$-bimodule linearity for $G$. Right linearity is similar. 
\end{proof}
\begin{remark} For $\mc{M}=\mc{C}$ the action (\ref{centralactiondef}) is the same as the monoidal structure for $\cent{C}{D}$ described in (\ref{relcentTens}).
\end{remark}
\begin{remark}\label{2-remark}
The $*$ composition in (\ref{centralactiondef}) can be described in a more general setting using 2-compositions as follows. Assume all cells occur in some fixed 2-category. Write $\alpha *_\sigma\beta$ for the 2-cell which corresponds to the 2-diagram below.
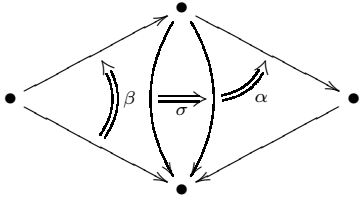
\[\xymatrix{
&&\bullet\ar[drr]_{}="4"\ar@/^1pc/[dd]^{}="3"_{}="6"\ar@/_1pc/[dd]^{}="5"&&\\
\bullet\ar[urr]_{}="2"\ar[drr]^{}="1"&&&&\bullet\ar[dll]\\
&&\bullet&&\\
\ar@/_.4pc/@{=>}"3";"4"_{\alpha} 
\ar@/_.6pc/@{=>}"1";"2"_{\beta} 
\ar@{=>}"5";"6"_{\sigma}}
\vspace{-2.5em}\]
This operation, which we will call \textit{convolution} of $\alpha$ and $\beta$ over $\sigma$, is associative: to show, for example, that $(\alpha*_\gamma\beta)*_\sigma\epsilon = \alpha*_\gamma(\beta*_\sigma\epsilon)$ one must simply check that the 2-cells corresponding to the diagrams
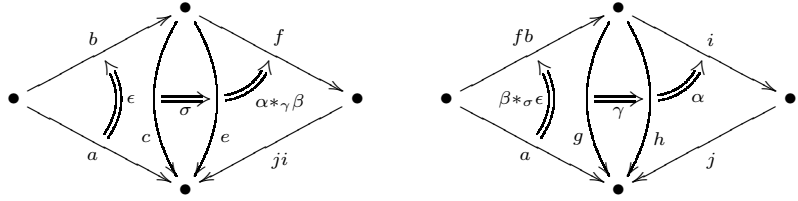
\[\xymatrix{
&&\bullet\ar[drr]_{}="4"^f\ar@/^1pc/[dd]^{}="3"_{}="6"^>>>>>>e\ar@/_1pc/[dd]^{}="5"_>>>>>>c&&\\
\bullet\ar[urr]_{}="2"^b\ar[drr]^{}="1"_a&&&&\bullet\ar[dll]^{ji}\\
&&\bullet&&\\
\ar@/_.4pc/@{=>}"3";"4"_{\alpha*_\gamma\beta} 
\ar@/_.6pc/@{=>}"1";"2"_{\epsilon} 
\ar@{=>}"5";"6"_{\sigma}}
\;\qquad
\xymatrix{
&&\bullet\ar[drr]_{}="4"^i\ar@/^1pc/[dd]^{}="3"_{}="6"^>>>>>>h\ar@/_1pc/[dd]^{}="5"_>>>>>>g&&\\
\bullet\ar[urr]_{}="2"^{fb}\ar[drr]^{}="1"_a&&&&\bullet\ar[dll]^j\\
&&\bullet&&\\
\ar@/_.4pc/@{=>}"3";"4"_{\alpha} 
\ar@/_.6pc/@{=>}"1";"2"^{\beta*_\sigma\epsilon} 
\ar@{=>}"5";"6"_{\gamma}}
\vspace{-2.5em}\]
agree. Using the relations available in any 2-category one shows without difficulty that they are both given by $((((\alpha\gamma)*f)\beta\sigma)*b)\epsilon$.

For $c:\otimes\stackrel{\sim}{\rightarrow}\otimes\tau:\mc{D}\times\mc{C}\rightarrow\mc{C}$ and $\varphi:\otimes\stackrel{\sim}{\rightarrow}\otimes\tau:\mc{C}\times\mc{M}\rightarrow\mc{M}$ the composition $c*\varphi$ may be represented as this convolution of 2-cells with $\beta=c\otimes id_\mc{M}$, $\alpha=id_\mc{C}\otimes\varphi$ and $\sigma = id$; $c*\varphi$ in (\ref{centralactiondef}) (really $(1\otimes\varphi)*_\alpha(c\otimes 1)$) is the 2-cell represented by the diagram
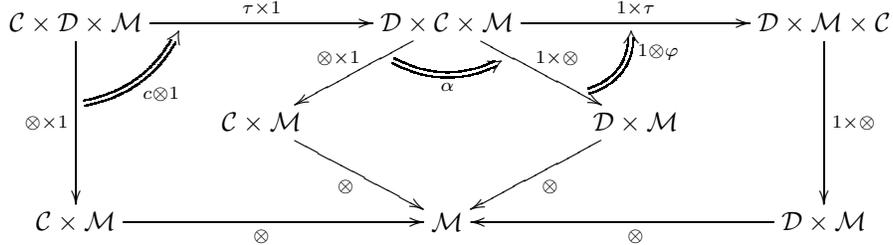
\[
\xymatrix{
\mc{C}\times\mc{D}\times\mc{M}\ar[rr]^{\tau\times 1}_<<<<<{}="2"\ar[dd]_{\otimes\times 1}^<<<<<<<<<{}="1"&&\mc{D}\times\mc{C}\times\mc{M}\ar[dl]_{\otimes\times 1}^<<<<<{}="3"\ar[dr]^{1\times\otimes}_<<<<<{}="4"^>>>>{}="5"\ar[rr]^{1\times\tau}_>>>>>>>>>>>>>>>>{}="6"&&\mc{D}\times\mc{M}\times\mc{C}\ar[dd]^{1\times\otimes}\\
&\mc{C}\times\mc{M}\ar[dr]_\otimes&&\mc{D}\times\mc{M}\ar[dl]^\otimes&\\
\mc{C}\times\mc{M}\ar[rr]_{\otimes}&&\mc{M}&&\mc{D}\times\mc{M}\ar[ll]^\otimes
\ar@/_.6pc/@{=>}"1";"2"_{c\otimes 1} 
\ar@/_.6pc/@{=>}"3";"4"_{\alpha} 
\ar@/_.6pc/@{=>}"5";"6"_<<<<<<<{1\otimes\varphi} 
}\]
The natural isomorphism $\alpha$ is trivial in case the $\mc{C}$ action is strict, and we assume that it is.
\end{remark}
\begin{cor}\label{firstDecomCor.1} Let $\mc{M}_i, \mc{N}_i, i=1, 2$ be module categories over tensor categories $\mc{C}_1, \mc{C}_2$, respectively, and suppose further that both $\mc{C}_i$ are tensor over $\mc{D}$. Further suppose that $F_i:\mc{M}_i\rightarrow\mc{N}_i$, $i=1, 2$, are a pair of functors where each $F_i$ is a $\mc{C}_i$-module functor. Then 
\begin{itemize}
\item[i)] $\bitens{M}{D}{N}$ has canonical structure of a $\mc{C}_1\tens{D}\mc{C}_2$-module category. 
\item[ii)] ${F_1}\tens{D}{F_2}:\mc{M}_1\tens{D}\mc{M}_2\rightarrow \mc{N}_1\tens{D}\mc{N}_2$ has canonical structure of a $\mc{C}_1\tens{D}\mc{C}_2$-module functor. 
\end{itemize}
In fact, $\tens{D}$ is a 2-functor $\mc{C}_1$-Mod $\tens{}\mc{C}_2$-Mod $\rightarrow\mc{C}_1\tens{D}\mc{C}_2$-Mod.
\end{cor}
\begin{thm}\label{2-functorCent} Let $\mc{C}$ be tensor over braided fusion category $\mc{D}$. Then the association $\mc{Z}_{\mc{D}}:\mc{C}$-Bimod $\rightarrow \cent{C}{D}$-Bimod is a 2-monoidal 2-functor. 
\end{thm}
\begin{proof} It remains to check that $\mc{Z}_\mc{D}$ operates in the appropriate way on 2-cells, and that the relevant composition rules are retained. We leave this straightforward but tedious exercise to the reader. 
\end{proof}
Proposition \ref{FunctorialityOfZ_DLemma} says that, under the stated hypotheses, $\cent{M}{D}$ has the structure of a $\cent{C}{D}$-module category with action described by equation (\ref{centralactiondef}). On the other hand $\cent{M}{C}$ has the natural structure of a left $\mc{D}$-module category with action
\begin{equation}
D\otimes (M, \varphi_M):=(D\otimes M, c_D*\varphi_M)
\end{equation}
where as above $(c_D*\varphi_M)_X:=(id_D\otimes\varphi_{X, M})(c_{X, D}\otimes id_M)$. Here $c_D$ denotes the braiding in $\mc{Z}(\mc{C})$. Strictness of both associativity of the action and the bimodule consistency of left and right $\mc{D}$-module actions (Definition \ref{BimoduleCatDef}) are implied by Remark \ref{2-remark}. We may therefore ask about the the relative center of $\cent{M}{C}$ over $\mc{D}$ and if this is at all related to the relative center of $\cent{M}{D}$ over $\cent{C}{D}$. We expect some description of how to intertwine the functors $\mc{Z}_\mc{D}$ and $\mc{Z}_{\cent{C}{D}}$ which reduces to \cite[Theorem 7.14]{JG} when $\mc{C}=\mc{D}$. The general result describing this interaction is slightly beyond our grasp at present, though the following proposition is a step in the direction of a solution via functoriality of the relative tensor product.
\begin{prop}\label{functorialityTensor}
Let $\mc{C}_i, i\in\{1, 2\}$ be tensor over $\mc{D}$. Suppose that $\mc{M}_i$ are right $\mc{C}_i$-module categories and $\mc{N}_i$ left $\mc{C}_i$-module categories. Then there is a canonical equivalence 
\begin{equation}\label{TensFunctoriality}
(\mc{M}_1\boxtimes_{\mc{C}_1}\mc{N}_1)\boxtimes_{\mc{D}}(\mc{M}_2\boxtimes_{\mc{C}_2}\mc{N}_2)\simeq(\mc{M}_1\tens{D}\mc{M}_2)\boxtimes_{\mc{C}_1\tens{D}\mc{C}_2}(\mc{N}_1\tens{D}\mc{N}_2).
\end{equation}
\end{prop}
\begin{proof} To prove the proposition we construct the diagram below.
\[\RQq\RQq\RQq\!\!\!
\xymatrix{
\mc{M}_1\tens{}\mc{N}_1\tens{}\mc{M}_2\tens{}\mc{N}_2\ar[r]^{(23)}\ar[d]_{ B_{\mc{M}_1, \mc{N}_1}\tens{} B_{\mc{M}_2, \mc{N}_2}}&
\mc{M}_1\tens{}\mc{M}_2\tens{}\mc{N}_1\tens{}\mc{N}_2\ar[r]^{(23)}\ar[d]_{B'_{\mc{M}}\tens{}B'_{\mc{N}}}&
\mc{M}_1\tens{}\mc{N}_1\tens{}\mc{M}_2\tens{}\mc{N}_2\ar[d]\\
%
%
\mc{M}_1\boxtimes_{\mc{C}_1}\mc{N}_1\tens{}\mc{M}_2\boxtimes_{\mc{C}_2}\mc{N}_2\ar[d]_{B}\ar@{-->}[dr]^{\overline{\mc{V}}}
&\mc{M}_1\tens{D}\mc{M}_2\tens{}\mc{N}_1\tens{D}\mc{N}_2\ar[d]_{B'}\ar@{-->}[dr]
&\mc{M}_1\boxtimes_{\mc{C}_1}\mc{N}_1\tens{}\mc{M}_2\boxtimes_{\mc{C}_2}\mc{N}_2\ar[d]\\
%
%
\mc{M}_1\boxtimes_{\mc{C}_1}\mc{N}_1\tens{D}\mc{M}_2\boxtimes_{\mc{C}_2}\mc{N}_2\ar[r]&
(\mc{M}_1\tens{D}\mc{M}_2)*(\mc{N}_1\tens{D}\mc{N}_2)\ar[r]&
\mc{M}_1\boxtimes_{\mc{C}_1}\mc{N}_1\tens{D}\mc{M}_2\boxtimes_{\mc{C}_2}\mc{N}_2
}\]
Here we have abbreviated $B :=B_{\mc{M}_1\boxtimes_{\mc{C}_1}\mc{M}_2, \mc{N}_1\boxtimes_{\mc{C}_2}\mc{N}_2}$, $B':=B_{\mc{M}_1\tens{D}\mc{M}_2, \mc{N}_1\tens{D}\mc{N}_2}$, and $*:=\boxtimes_{\mc{C}_1\tens{D}\mc{C}_2}$ in the interest of space. The third column is identical to the first. 

We begin by showing that the composition of functors $\mc{V}:=B(B'_{\mc{M}}\tens{}B'_{\mc{N}})(23)$ in the diagram is $\mc{C}_1$-balanced with respect to the first $\tens{}$ and is $\mc{C}_2$-balanced with respect to the third. Let $M_i\in\mc{M}_i$, $N_i\in\mc{N}_i$ and $X\in\mc{C}_1$, and abbreviate the image of the universal functors $B$ by tensor product: $B_{\mc{M}, \mc{N}}(M\tens{}N)=M\tens{D}N$. Then
\begin{eqnarray*}
\mc{V}((M_1\otimes X)\tens{}N_1\tens{}M_2\tens{}N_1)&=&B((M_1\tens{D}M_2)\otimes (X\tens{D}1)\tens{}(N_1\tens{D}N_2))\\
&\simeq&B((M_1\tens{D}M_2)\tens{}(X\tens{D}1)\otimes(N_1\tens{D}N_2))\\
&=&\mc{V}(M_1\tens{}(X\otimes N_1)\tens{}M_2\tens{}N_1).
\end{eqnarray*}
The first equation is definition of $\mc{V}$ and the $\mc{C}_1\tens{D}\mc{C}_2$-module structure in $\mc{M}_1\tens{D}\mc{M}_2$. The equivalence in the second line is the $\mc{C}_1\tens{D}\mc{C}_2$-balancing of the universal functor $B$ occurring in the middle column of the diagram. One similarly shows that $\mc{V}$ is $\mc{C}_2$-balanced in the third position, giving a unique functor $\overline{\mc{V}}$ making the pentagonal subdiagram commute. 

Next we show that $\overline{\mc{V}}$ is $\mc{D}$-balanced. Give the module categories $\mc{M}_i$, $\mc{N}_i$ the structure of $\mc{D}$-bimodule categories using $\mc{C}_i$-module category structures and the fact that $\mc{D}$ is braided. The right $\mc{D}$-module action on $\mc{M}\boxtimes_{\mc{C}_1}\mc{N}_1$ is given by $(M\boxtimes_{\mc{C}_1}N)\otimes D=M\boxtimes_{\mc{C}_1}(N\otimes D)$ on ``simple tensors." Thus
\begin{eqnarray*}\!\!\!\!\!\!
\overline{\mc{V}}(((M_1\boxtimes_{\mc{C}_1}N_1)\otimes D)\tens{}(M_2\boxtimes_{\mc{C}_2}N_2))&=&B(B'_{\mc{M}}(M_1\tens{}M_2)\tens{}B'_{\mc{N}}((N_1\otimes D)\tens{}N_2))\\
&\simeq&B(B'_{\mc{M}}(M_1\tens{}M_2)\tens{}B'_{\mc{N}}(N_1\tens{}(D\otimes N_2)))\\
&=&\overline{\mc{V}}((M_1\boxtimes_{\mc{C}_1}N_1)\tens{}(M_2\boxtimes_{\mc{C}_2}(D\otimes N_2)))\\
&\simeq&\overline{\mc{V}}((M_1\boxtimes_{\mc{C}_1}N_1)\tens{}((D\otimes M_2)\boxtimes_{\mc{C}_2}N_2)).
\end{eqnarray*}
First and third lines are the definition of $\overline{\mc{V}}$, second is the $\mc{D}$-balancing of $B'_{\mc{N}}$ and the fourth follows from the $\mc{C}_2$-balancing of $B_{\mc{M}_2, \mc{N}_2}$ and the $\mc{D}$-bimodule structure in $\mc{M}_2$. Hence $\overline{\mc{V}}$ is a $\mc{D}$-balanced functor. We therefore get a unique functor (unlabeled horizontal arrow in lower left) making the rectangular left half of the large diagram commute. 

Similarly, one shows that $B_{(1, 2)}(B_{\mc{M}_1, \mc{N}_1}\tens{} B_{\mc{M}_2, \mc{N}_2})(23)$ is $\mc{D}$-balanced in the first and third position using the balanced structures in $B_{\mc{M}_i, \mc{N}_i}$ and $B_{(1, 2)}$. This yields a unique functor (unlabeled perforated arrow in right side of diagram) making the resulting pentagonal subdiagram commute.  Using balancing of $B'_\mc{M}$, $B'_\mc{N}$ and the definition of the $\mc{C}_1\tens{D}\mc{C}_2$-module category action in $\mc{M}_1\tens{D}\mc{M}_2$ and $\mc{N}_1\tens{D}\mc{N}_2$ it's easy to show that this functor is itself $\mc{C}_1\tens{D}\mc{C}_2$-balanced, giving rise to a unique functor (the unlabeled solid arrow in lower right) making the resulting triangle commute. On uniqueness and $(23)^2=id$ the composition along the bottom of the diagram is identity, giving the equivalence appearing in the statement of the lemma.
\end{proof}
\begin{note}
In the case that categories $\mc{M}_i, \mc{N}_i$ have bimodule category structure the equivalence in Proposition \ref{functorialityTensor} is an equivalence of bimodule categories. 
\end{note}
%
%
\begin{cor}\label{tensorinversion} Suppose that $\mc{M}_{1}$ is $\mc{C}_1$-bimodule category, $\mc{M}_{2}$ is a $\mc{C}_2$-bimodule category, and that the $\mc{C}_i$ are both tensor over $\mc{D}$. Then if $\mc{M}_{i}$ are invertible as $\mc{C}_i$-bimodule categories then $\mc{M}_1\tens{D}\mc{M}_2$ is invertible as a ${\mc{C}_1\tens{D}\mc{C}_2}$-bimodule category and has inverse $(\mc{M}_1\tens{D}\mc{M}_2)^{-1}=\mc{M}_{1}^{-1}\tens{D}\mc{M}_{2}^{-1}$.
\end{cor}
\begin{proof} A straightforward application of Proposition \ref{functorialityTensor} gives
\begin{equation*}
(\mc{M}_{1}\boxtimes_{\mc{D}}\mc{M}_{2})\boxtimes_{{\mc{C}_1\tens{D}\mc{C}_2}}(\mc{M}_{1}^{-1}\boxtimes_{\mc{D}}\mc{M}_{2}^{-1})\simeq(\mc{M}_{1}\boxtimes_{\mc{C}_1}\mc{M}_{1}^{-1})\tens{D}(\mc{M}_{2}\boxtimes_{\mc{C}_2}\mc{M}_{2}^{-1})
\end{equation*}
and the second term is equivalent to $\mc{C}_1\tens{D}\mc{C}_2$ by definition.
\end{proof}
\section{Centers and tensor product over centralizers}\label{AlmostNonDegSect}
Let $\mc{C}$ be a tensor category, and let  $\mc{C}_1, \mc{C}_2$ be tensor subcategories of tensor category $\mc{C}$. Following  \cite{braidedI} we denote by $\mc{C}_1\vee\mc{C}_2$ the smallest tensor subcategory of $\mc{C}$ containing both $\mc{C}_1$ and $\mc{C}_2$. Recall that if $\mc{C}$ is braided with braiding $c$ then $\mc{C}'$ denotes the braided tensor subcategory of $\mc{C}$ consisting of objects $X$ such that
\begin{equation}\label{someEquation}
c_{Y, X}c_{X, Y}=id_{X\otimes Y}
\end{equation}
for all $Y\in\mc{C}$. We call this tensor subcategory the \textit{centralizer} of $\mc{C}$. 

Now suppose $\mc{C}$ is braided. There are two natural braided inclusions 
\begin{eqnarray}
\eta_+:\mc{C}\hookrightarrow \mc{Z}(\mc{C}),&&\quad X\mapsto (X, c_{\score, X})\label{eta1}\\
\eta_-:\mc{C}^{rev}\hookrightarrow \mc{Z}(\mc{C}),&&\quad X\mapsto (X, c^{-1}_{X, \score}).\label{eta2}
\end{eqnarray}
$\mc{C}^{rev}$ is the tensor category which is $\mc{C}$ as an abelian category but with opposite tensor structure. The images of $\eta_\pm$ form braided tensor subcategories of $\mc{Z}(\mc{C})$ which we abbreviate $\mc{C}_\pm$, respectively. Verification that $\mc{C}_\pm$ are $\mc{C}'$-bimodule categories is straightforward and left to the reader. We note here that the left and right $\mc{C}'$-module actions agree, a fact which relies on (\ref{someEquation}).
\begin{lem}\label{minimality of tensor product} Let $\mc{C}$, $\mc{D}$ be braided subcategories of a braided category $\mc{B}$ such that both $\mc{C}$ and $\mc{D}$ are braided over $\mc{C}\cap\mc{D}$. Then $\mc{C}\vee\mc{D}\simeq\mc{C}\boxtimes_{\mc{C}\cap\mc{D}}\mc{D}$ as braided tensor categories. If all categories in sight are fusion this is an equivalence of braided fusion categories.
\end{lem}
\begin{proof} Note first that $\mc{C}$, $\mc{D}$ being braided over $\mc{C}\cap\mc{D}$ is equivalent to $\mc{C}\cap\mc{D}\subset\mc{C}'\cap\mc{D}'$, the latter being symmetric. Hence the braiding on $\mc{C}\boxtimes_{\mc{C}\cap\mc{D}}\mc{D}$ is well defined (see Definition \ref{braidingOVERdef} and Theorem \ref{braiding of relative center thm}). 

To say that $\mc{C}\vee\mc{D}$ is the \textit{smallest} braided tensor subcategory of $\mc{B}$ containing both $\mc{C}$ and $\mc{D}$ as tensor subcategories is to say that, given any pair of braided inclusions $v_1, v_2$ into some other braided subcategory $\mc{V}$ of $\mc{B}$, there is a unique braided inclusion $Q$ making the diagram below commute. 
\[
\xymatrix{
\mc{C}\ar@{->}[r]^{v_1}\ar[dr]&\mc{V}&\mc{D}\ar[l]_{v_2}\ar[dl]\\
&\mc{C}\vee\mc{D}\ar[u]_Q&
}
\]
The diagonal arrows are the canonical braided inclusions. We show that the relative tensor product $\mc{C}\boxtimes_{\mc{C}\cap\mc{D}}\mc{D}$ satisfies the unique minimality property defining $\mc{C}\vee\mc{D}$. Define braided inclusions 
\begin{equation*}
\eta_1:\mc{C}\rightarrow\mc{C}\tens{}\mc{D},\,\, X\mapsto X\tens{}1\qquad\eta_2:\mc{D}\rightarrow\mc{C}\tens{}\mc{D},\,\, Y\mapsto 1\tens{}Y.
\end{equation*}
Then the $\eta_i$ and the braided inclusions $v_i$ combine to give braided inclusions 
\begin{equation*}
\otimes_\mc{V}:\mc{C}\tens{}\mc{D}\rightarrow\mc{V},\quad X\tens{}Y\mapsto v_1(X)\otimes v_2(Y).
\end{equation*}
Since $v_1$ and $v_2$ agree on $\mc{C}\cap\mc{D}$ the functor $\otimes_\mc{V}$ is $\mc{C}\cap\mc{D}$-balanced, and therefore descends to a unique braided inclusion $\overline{\otimes}_\mc{V}:\mc{C}\boxtimes_{\mc{C}\cap\mc{D}}\mc{D}\rightarrow\mc{V}$. Everything in this discussion so far is contained in the commuting diagram below. 
\[
\xymatrix{
\mc{C}\ar[rr]^{v_1}\ar[d]_{\eta_1}&&\mc{V}&&\mc{D}\ar[ll]_{v_2}\ar[d]^{\eta_2}\\
\bitens{C}{}{D}\ar[urr]^{\otimes_\mc{V}}\ar[rr]_{B}&&\mc{C}\boxtimes_{\mc{C}\cap\mc{D}}\mc{D}\ar[u]_{\overline{\otimes}_\mc{V}}&&\bitens{C}{}{D}\ar[ll]^{B}\ar[ull]_{\otimes_\mc{V}}
}\]
The functor $B$ is the universal $\mc{C}\cap\mc{D}$-balanced functor. It is easy to check that the compositions $B\eta_i$ are inclusions and thus by Theorem \ref{braiding of relative center thm} are braided inclusions. On the minimality of the category $\mc{C}\vee\mc{D}$ we obtain the equivalence stated in the lemma. 
\end{proof}
\begin{prop}\label{tensorveeprop}  $\mc{C}_+\vee\mc{C}_-\simeq\mc{C}\tens{C'}\mc{C}^{rev}$ canonically as braided tensor categories.
\end{prop}
\begin{proof}
This is just Lemma \ref{minimality of tensor product} with $\eta_{\pm}$ replacing $\eta_1, \eta_2$, and  observing that $\mc{C}\tens{}1=\mc{C}_+\simeq\mc{C}$, $1\tens{}\mc{C}^{rev}=\mc{C}_-\simeq\mc{C}^{rev}$ and that $\mc{C}'=\mc{C}_+\cap\mc{C}_-$.
\end{proof}
By Proposition \ref{tensorveeprop} we may identify $\widetilde{\mc{C}}:=\mc{C}\boxtimes_{\mc{C}'}\mc{C}^{rev}$ with a braided fusion subcategory of $\mc{Z}(\mc{C})$. As a $\widetilde{\mc{C}}$-bimodule category we therefore have the decomposition
\begin{equation}\label{CenterDecompEquation}
\mc{Z}(\mc{C})\simeq\widetilde{\mc{C}}\oplus\mc{R}
\end{equation}
for some $\widetilde{\mc{C}}$-bimodule subcategory $\mc{R}$ of $\mc{Z}(\mc{C})$. A complete description of the center $\mc{Z}(\mc{C})$ as a $\widetilde{\mc{C}}$-bimodule category amounts to a description of this nontrivial component $\mc{R}$. 
\begin{lem}\label{tensorveeFPlem} With notation as above $\FPdim(\mc{Z}(\mc{C}))=\FPdim(\widetilde{\mc{C}})\FPdim(\mc{C}')$.
\end{lem}
\begin{proof} Lemma 3.38 in \textit{loc. cit.} and Proposition \ref{tensorveeprop} imply 
\begin{equation*}
\FPdim(\widetilde{\mc{C}})=\frac{\FPdim(\mc{C}_+)\FPdim(\mc{C}_-)}{\FPdim(\mc{C}_+\cap\mc{C}_-)}.
\end{equation*}
The result follows by noting that $\mc{C}_+\cap\mc{C}_-=\mc{C}'$, $\FPdim(\mc{C}_\pm)=\FPdim(\mc{C})$ and recalling that $\FPdim(\mc{Z}(\mc{C}))=\FPdim(\mc{C})^2$.
\end{proof}
\begin{ex}\label{groupGradedExample} If the center $\mc{Z}(\mc{C})$ is faithfully graded by a group $G$ as a $\mc{C}_+\vee\mc{C}_-$-bimodule category then $\FPdim(\mc{C}')=|G|$. In particular, if $\mc{C}$ is non-degenerate then $\mc{Z}(\mc{C})=\mc{C}_+\vee\mc{C}_-$. Indeed, applying $\FPdim$ to both sides of (\ref{CenterDecompEquation}) and using Lemma \ref{tensorveeFPlem} gives
\begin{equation}\label{FPdimFormula}
\FPdim(\mc{C})^2(1-(\FPdim(\mc{C}')^{-1})=\FPdim(\mc{R}).
\end{equation}
In the case $\mc{Z}(\mc{C})=\oplus_{g\in G}\mc{C}_g$ for $\mc{C}_e=\mc{C}_+\vee\mc{C}_-$ Proposition 8.20 in \cite{ENO:ofc} implies that $\FPdim(\mc{C}_g)=\FPdim(\mc{C})^2/\FPdim(\mc{C}')$ for every element $g\in G$. Solving for $\FPdim(\mc{C}')$ in (\ref{FPdimFormula}) using $\FPdim(\mc{R})=(|G|-1)\FPdim(\mc{C}_e)$ gives the result. In the special case that $\mc{C}'=Vec$ the center is graded by the trivial group, i.e. has a single irreducible component as a bimodule category.
\end{ex}
\subsection{An in-depth example: $\mc{Z}(\Rep(G))$ as a $\Rep(G)$-bimodule category.} In this section we examine the decomposition (\ref{CenterDecompEquation}) for the case $\mc{C}=\Rep(G)$, $G$ a finite group. It is known that $(Vec_G)^*_{Vec}=\Rep(G)$ and hence $\mc{Z}(\Rep(G))\simeq\mc{Z}(Vec_G)$. Because graded vector spaces are in some senses easier to work with we will discuss the center of $\Rep(G)$ in terms of $Vec_G$.

It is known that the category $\mc{Z}(Vec_G)$ is equivalent to the category of $G$-equivariant vector bundles over the group $G$ with respect to the action of $G$ on itself by conjugation (see e.g. \cite{springerlink:10.1007/BF02365308}). That is, objects of $\mc{Z}(Vec_G)$ are pairs $(V, \gamma)$ where $V$ is a $G$-graded vector space and $\gamma_g:g\otimes V\rightarrow V\otimes g$ are natural isomorphisms (we have abbreviated simple objects in $Vec_G$ by group elements). The natural isomorphisms $\gamma$ satisfy the usual pair of hexagons appearing in the definition of a braided category. On the level of components at summand $s\in G$ these reduce to
\begin{equation}\label{someOtherEquation}
(\gamma_g\otimes h)_s(g\otimes\gamma_h)_s = (\gamma_{gh})_s
\end{equation}
for any pair $g, h\in G$ since $Vec_G$ is strict. Now define linear isomorphisms
\begin{equation}\label{piComponentIsos}
\pi_{g, x} : = (\gamma_{g^{-1}})_{xg}:V_{g^{-1}xg}\rightarrow V_x.
\end{equation}
Explicitly, $\pi_{g, x}$ is the $xg^{th}$ component of $\gamma_{g^{-1}}$. Rewriting (\ref{someOtherEquation}) in terms of $\pi$ gives the relation $\pi_{gh, x} = \pi_{g, x}\circ\pi_{h, g^{-1}xg}$. The isomorphisms (\ref{piComponentIsos}) imply, in particular, that graded components of objects in $\mc{Z}(Vec_G)$ are (up to isomorphism) constant on conjugacy classes. 

Consider conjugacy class $\overline{a}$ containing $a\in G$ and denote by $Z(a)$ the centralizer of $a$ in $G$. Then for $g, h\in Z(a)$ $\pi_{gh, a}=\pi_{g, a}\circ\pi_{h, a}$ and the assignment $\rho:Z(a)\rightarrow\End (V_a)$, $g\mapsto\pi_{g, a}$ gives an irreducible representation of $Z(a)$. In this way simple objects of $\mc{Z}(Vec_G)$ correspond to pairs $(\overline{a}, \rho)$ where $\rho$ is an irreducible representation of the centralizer $Z(a)$ of $a$ (here $\pi_{\underline{\:\:}, a} = \rho$). As a graded object of the center $(\overline{a}, \rho)$ is the object which has a copy of $\rho$ at the $x$th place for every $x\in\overline{a}$ and $0$ elsewhere. 

Denote by $\mc{C}_a$ the subcategory of $\mc{Z}(Vec_G)$ consisting of objects supported on $\overline{a}$, and let $\mc{O}_a$ be the set of irreducible objects of $\Rep(Z(a))$. The induction functor $\Rep(Z(a))\rightarrow\mc{C}_a$, $V\mapsto V^G$, constitutes a bijection $\mc{O}_a\rightarrow\mc{O}(\mc{C}_a)$: for $V\in\mc{O}_a$ one has $V^G=(\overline{a}, V)$ since $|G/Z(a)|=|\overline{a}|$, and therefore $\FPdim(V^G)=|\overline{a}|\dim(V)$.

The category $\Rep(G)$ sits inside $\mc{Z}(Vec_G)$ (non uniquely) as the braided subcategory $\mc{C}_q$ for any $q\in Z(G)$; for any irreducible representation $\rho$ of $G=Z(q)$ send $\rho\mapsto (q, \rho)$. The monoidal structure in $\mc{Z}(Vec_G)$ then determines a $\Rep(G)$-module category structure. Since $\Rep(G)$ is symmetric we have $\Rep(G)'=\Rep(G)$, hence $\widetilde{\Rep(G)}\simeq\Rep(G)$ as fusion categories. The decomposition of $\mc{Z}(Vec_G)$ into indecomposable $\Rep(G)$-module categories will in fact be the decomposition appearing in (\ref{CenterDecompEquation}). We have the following lemma.
\begin{lem}
For $H<G$ denote by $\Rep(H)^G$ the subcategory of $\mc{Z}(Vec_G)$ generated by simple objects of the form $V^G$ for irreducible representations $V$ of $H$. Then $\Rep(Z(a))^G$ and $\mc{C}_a$ are equivalent as $\Rep(G)$-module categories. 
\end{lem}
Furthermore, since irreducible $\Rep(G)$-module categories correspond to categories of projective representations of subgroups of $G$ (\cite[Section 8]{JG}), the categories $\mc{C}_a$ are irreducible as $\Rep(G)$-module categories. 

The next step is to decompose $\mc{Z}(\Rep(G))$ into a sum of these $\mc{C}_a$. Let $R_G$ be the regular object of $Vec_G$ and define $R_a :=\mc{C}_a\cap R_G$, i.e. $R_a = \oplus_{V\in\mc{O}_a}\dim(V^G)V^G$. Then taking the Frobenius-Perron dimension of both sides we obtain $\FPdim(\mc{C}_a) = \sum_{V\in\mc{O}_a}|\overline{a}|^2\dim(V)^2$. Summing over all $\FPdim(\mc{C}_a)$ allowing each to occur exactly once for every $a\in R:=\{$conjugacy classes of $G\}$ gives 
\begin{equation*}
\sum_{r\in R}|\overline{r}|^2\sum_{V\in\mc{O}_r}\dim(V)^2 =\sum_{r\in R}|\overline{r}|^2|Z(r)| = |G|^2 = \FPdim\,\mc{Z}(\Rep(G)).
\end{equation*}
Since $\FPdim(\mc{C}_r)\neq 0$ for every $r$ we have, as $\Rep(G)$-module categories,
\begin{equation}\label{central decomp}
\mc{Z}(\Rep(G)) \simeq\bigoplus_{r\in R}\mc{C}_r
\end{equation}
where $R$ is taken to be a complete set of representatives of conjugacy classes in $G$. We record this observation in the following theorem. 
\begin{thm}\label{ClassIndexedIndecomposibles}
As a $\widetilde{\Rep(G)}$-module category $\mc{Z}(\Rep(G))$ decomposes into indecomposable module subcategories indexed by the conjugacy classes of $G$.
\end{thm}
\begin{remark}Theorem \ref{ClassIndexedIndecomposibles} follows also by Corollaries 3.6, 3.9 in \cite{braidedI}. 
\end{remark}
Now since, for any $q\in Z(G)$, the subcategory $\mc{C}_q\simeq\Rep(G)$, Theorem \ref{ClassIndexedIndecomposibles} allows us to classify the non-trivial component $\mc{R}$ in(\ref{CenterDecompEquation}) up to equivalence: 
\begin{equation*}
\mc{R}\simeq\bigoplus_{r\in R^*}\mc{C}_r
\end{equation*}
where $R^*$ is the set of representatives of conjugacy classes in $G$ minus that containing the identity. 

The next proposition tells us how to multiply the irreducible components appearing in the decomposition (\ref{central decomp}). 
\begin{prop}[$\mc{Z}(\Rep(G))$ fusion rules]\label{ConjugSumsStructProp} As $\Rep(G)$-module categories
\begin{equation*}
\mc{C}_a\boxtimes_{\Rep(G)}\mc{C}_b = \bigoplus_{r\in R}N^{ab}_r\mc{C}_r
\end{equation*}
where $N^{ab}_s$ are the structure coefficients in the ring of conjugacy sums over the group ring $\mathbb{Z}[G]$, i.e. $N^{ab}_s =|\{(x, y)\in \overline{a}\times\overline{b}|xy\in\overline{s}\}|$.
\end{prop}
\begin{note} Proposition \ref{ConjugSumsStructProp} states that the ring of isomorphism classes of indecomposable $\Rep(G)$-module subcategories of $\mc{Z}(\Rep(G))$ has structure coefficients which given by those defining the structure of the center of the group ring $\mathbb{Z}[G]$.
\end{note}
\begin{proof}
We begin by defining a functor $\mc{C}_a\tens{}\mc{C}_b\rightarrow\Rep(G)$ using the tensor product in $\Rep(G)$: $(\overline{a}, V)\tens{}(\overline{b}, W)\mapsto\ind^G_{Z(a)}(V)\otimes\ind^G_{Z(b)}(W)$. Rewriting each of these induced representations as their associated direct sums their tensor product becomes 
\begin{equation}\label{CenterDecompDoubles}
\bigoplus_{x\in\overline{a}}xV\otimes\bigoplus_{y\in\overline{b}}y W=\bigoplus_{g\in G}\bigoplus_{{xy=g}\atop{{x\sim a}\atop{y\sim b}}}(xV\otimes yW).
\end{equation}
Write $xV\otimes yW = xy(V\otimes W)$. Then if $g\sim h$ in $G$ it is easy to see that $g(V\otimes W) = h(V\otimes W)$ up to isomorphism. Thus we can group the summands on the right in (\ref{CenterDecompDoubles}) according to conjugacy classes in $G$ giving the first equation in (\ref{FinalEquation}) below. 
\begin{equation}\label{FinalEquation}
\bigoplus_{g\in G}\bigoplus_{{xy=g}\atop{{x\sim a}\atop{y\sim b}}}(xV\otimes yW) = \bigoplus_{r\in R}\left(\bigoplus_{g\sim r}g(V\otimes W)\right)\sum_{{xy\sim r}\atop{{x\sim a}\atop{y\sim b}}}1=\bigoplus_{r\in R}N^{ab}_r(\overline{r}, V\otimes W)
\end{equation}
The direct sum indexed by terms congruent to a fixed $r\in R$ is isomorphic to $\ind_{Z(r)}^G(V\otimes W)$ as a representation of $G$. The sum indexed by the stack of congruences is precisely the coefficient of the basis element $e_{r}$ in $Z(\mathbb{Z}[G])$ occurring in the product $e_ae_b$ of the basis elements given by conjugacy sums over $a, b$ in $G$, here denoted $N^{ab}_r$. This explains the second equation in (\ref{FinalEquation}). 

Hence the tensor product takes $\mc{C}_a\tens{}\mc{C}_b\mapsto\oplus_{r\in R}N^{ab}_r\mc{C}_r$. This functor is evidently $\Rep(G)$-balanced (use associativity of the monoidal structure in $\mc{Z}(\Rep(G))$ or else directly using the relations satisfied by the coefficients $N^{ij}_k$ expressing associativity of the multiplication in $Z(\mathbb{Z}[G])$) and thus induces a functor $\mc{C}_a\boxtimes_{\Rep(G)}\mc{C}_b\mapsto\oplus_{r\in R}N^{ab}_r\mc{C}_r$ which is an equivalence because the $\mc{C}_r$ are indecomposable and because the coefficients $N^{ab}_r$ are explicitly obtained by evaluation on simple objects.
\end{proof}
In a subsequent article we will define and study a particular based ring generated by the terms in the partition of the regular virtual object of $\mc{Z}(\Rep(G))$ coming from the decomposition in Equation (\ref{central decomp}). This will allow us to interpret classical formulas from the theory of representation of finite groups in category-theoretical terms. 


\end{document}